\documentclass[11pt]{amsart}
\usepackage[tmargin=1in,lmargin=1in,bmargin=1in,rmargin=1in]{geometry}
\usepackage{amsfonts, amsmath, amsthm, amssymb}
\usepackage{verbatim}
\usepackage[outdir=./]{epstopdf}

\usepackage{graphicx}
\usepackage[hidelinks]{hyperref}
\newcommand{\R}{\mathbf{R}}

\newtheorem{hyp}{Hypothesis}

\newtheorem{theorem}{Theorem}[section]
\newtheorem{proposition}[theorem]{Proposition}
\newtheorem{lemma}[theorem]{Lemma}
\newtheorem{corollary}[theorem]{Corollary}

\newtheorem*{remark}{Remark}

\newcommand{\supp}{\operatorname{supp}}

\newif\ifdraft
\draftfalse

\begin{document}

\ifdraft
\title[Mass-critical refined Strichartz (Draft)]{Mass-critical refined Strichartz estimates 
for Schr\"{o}dinger operators}
\else
\title[Mass-critical refined Strichartz]{Refined mass-critical Strichartz estimates 
  for Schr\"{o}dinger operators}
\fi
\address{Department of
  Mathematics, UC Berkeley, Berkeley, CA 94720, USA.}
\author{Casey Jao}

\begin{abstract}
  We develop refined Strichartz estimates at $L^2$ regularity for a
  class of time-dependent Schr\"{o}dinger operators. Such refinements
  quantify near-optimizers of the Strichartz
  estimate, and play a pivotal part in the global theory of
  mass-critical NLS.  On one hand, the harmonic analysis is quite
  subtle in the $L^2$-critical setting due to an enormous group of
  symmetries, while on the other hand, the spacetime Fourier analysis
  employed by the existing approaches to the constant-coefficient
  equation are not adapted to nontranslation-invariant situations,
  especially with potentials as large as those considered in this
  article.

  Using phase space techniques, we reduce to proving certain analogues
  of (adjoint) bilinear Fourier restriction estimates. Then we extend
  Tao's bilinear restriction estimate for paraboloids to more general
  Schr\"{o}dinger operators. As a particular application, the
  resulting inverse Strichartz theorem and profile decompositions
  constitute a key harmonic analysis input for studying large data
  solutions to the $L^2$-critical NLS with a harmonic oscillator
  potential in dimensions $\ge 2$. This article builds on recent work
  of Killip, Visan, and the author in one space dimension.

\end{abstract}

\maketitle

\section{Introduction}


In this article, we prove sharpened forms of the Strichartz inequality
for nontranslation-invariant linear Schr\"{o}dinger equations with
$L^2$ initial data.  Recall that solutions to the linear
constant-coefficient Schr\"{o}dinger equation
\begin{align}
  \label{e:intro_free_particle}
  i\partial_t u = - \frac{1}{2}\Delta u, \quad u(0, \cdot ) = u_0 \in 
  L^2(\mathbf{R}^d),
\end{align}
satisfy the Strichartz inequality~\cite{Strichartz1977}
\begin{align}
  \label{e:intro_free_particle_str}
  \|u\|_{L^{\frac{2(d+2)}{d}}_{t,x}(\mathbf{R} \times \mathbf{R}^d)} \le C \|
  u(0, \cdot) \|_{L^2(\mathbf{R}^d)}.
\end{align}
On the other hand, it is also known if $u$ a solution that comes close
to saturating this inequality, then it must exhibit some
``concentration'';
see~\cite{carles-keraani,merle-vega,Moyua1999,begout-vargas}. Such
inverse theorems may be equivalently formulated as a refined
estimate
\begin{align}
  \label{e:intro_refined_str}
  \|u\|_{L^{\frac{2(d+2)}{d}}_{t,x}} \lesssim \| u\|_{X}^{\theta} \|u(0,
  \cdot)\|_{L^2(\mathbf{R}^d)}^{1-\theta},
\end{align}
where the norm $X$ is weaker than the right side
of~\eqref{e:intro_free_particle_str} but measures the ``microlocal
concentration'' of the solution.  We pursue analogues of such
refinements when the right side of~\eqref{e:intro_free_particle} is
replaced by a more general Schr\"{o}dinger operator
$-\tfrac{1}{2}\Delta + V(t,x)$.

Inverse theorems for the Strichartz inequality have provided a key
input to the study of the $L^2$-critical NLS
\begin{align}
  \label{e:intro_mass-crit_nls}
  i\partial_t u = -\frac{1}{2}\Delta u \pm |u|^{\frac{4}{d}} u, \quad u(0,
  \cdot) \in L^2(\mathbf{R}^d),
\end{align}
so termed because the rescaling
$u \mapsto u_{\lambda}(t,x) := \lambda^{d/2} u(\lambda^2 t, \lambda
x)$ preserves both the equation~\eqref{e:intro_free_particle} and the
$L^2$-norm
$M[u] := \| u(t)\|_{L^2(\mathbf{R}^d)} =
\|u(0)\|_{L^2(\mathbf{R}^d)}$. Indeed, they are used to construct the
profile decompositions underpinning the Bourgain-Kenig-Merle
concentration compactness and rigidity method by identifying potential
blowup scenarios for nonlinear solutions with large data. Using this
method, the large data global regularity problem
for~\eqref{e:intro_mass-crit_nls} was recently settled by
Dodson~\cite{dodson_mass-crit_1d,dodson_mass-crit_2d,dodson_mass-crit_high-d,dodson_focusing},
building on earlier work of Killip, Visan, Tao, and
Zhang~\cite{Killip2009, Killip2008, Tao2007}. For further discussion
of this equation we refer the interested reader to the lecture
notes~\cite{claynotes}.

The large group of symmetries for the
inequality~\eqref{e:intro_free_particle_str} is a significant obstruction to
characterizing its near-optimizers. Besides translation and scaling
symmetry, both sides are also invariant under \emph{Galilei}
transformations
\begin{align*}
  u \mapsto u_{\xi_0}(t,x) := e^{ i [\langle x, \xi_0 \rangle -
  \frac{1}{2}t|\xi_0|^2]} u (t, x-t\xi_0), \quad \xi_0 \in \mathbf{R}^d.
\end{align*}
%
This last symmetry emerges only at $L^2$ regularity and creates an
additional layer of complexity. In particular, while the Littlewood-Paley
decomposition is extremely well-adapted to higher Sobolev regularity variants
of~\eqref{e:intro_free_particle_str}, such as the $\dot{H}^1$-critical
estimate
\begin{align*}
  \|u\|_{L^{\frac{2(d+2)}{d-2}}_{t,x}} \lesssim \| \nabla 
  u(0)\|_{L^2(\mathbf{R}^d)},
\end{align*}
it is useless for inverting the $L^2$-critical estimate because one
has no a priori knowledge of where the solution is concentrated in
frequency. Instead, the mass-critical refinements
cited above combine spacetime Fourier-analytic arguments with
restriction theory for the paraboloid.

In physical applications, one is naturally led to consider variants of the mass-critical
equation~\eqref{e:intro_mass-crit_nls} with external potentials, such
as the harmonic oscillator
\begin{align}
  \label{e:intro_mass-crit_V}
  i\partial_t u = \Bigl( -\frac{1}{2}\Delta + \sum_j \omega_j^2 x_j^2 \Bigr) u 
  \pm |u|^{\frac{4}{d}} u, \quad u(0, \cdot) \in L^2(\mathbf{R}^d).
\end{align}
For instance, the cubic equation (with a $|u|^2u$ nonlinearity) has
been proposed as a model for Bose-Einstein condensates in a laboratory
trap~\cite{zhang_bec} where $\|u(t)\|_{L^2}$ represents the total 
number of particles, and in two space dimensions the critical
Sobolev norm for this equation is precisely $L^2$.

While introducing the potential breaks scaling symmetry, one
nonetheless expects solutions with highly concentrated
initial data to be approximated, for short times, by solutions to the
scale-invariant equation~\eqref{e:intro_mass-crit_nls}. Less
obviously, the equation is invariant under ``generalized'' Galilei
boosts, detailed in Lemma~\ref{l:galilei} below, where the
spatial and frequency parameters act together on the solutions; in the constant 
coefficient setting, this
reduces to the usual independent space translation and Galilei boost
symmetries.

This article 
develops
refined Strichartz estimates for the linear equation
\begin{align*}
  i\partial_t u = \Bigl( -\frac{1}{2}\Delta + V \Bigr) u, \quad u(0,
  \cdot) \in L^2(\mathbf{R}^d),
\end{align*}
for a class $\mathcal{V}$ of real-valued potentials $V(t, x)$ that
merely satisfy similar bounds as the harmonic oscillator and possibly
also depend on time.  Specifically, define
\begin{align}
  \label{e:potential_class_def}
  \mathcal{V} := \{ V : \R \times \R^d \to \R: \|
  \partial^\alpha_x V\|_{L^\infty_{t,x}} \le M_{|\alpha|} \text{ for } 2
  \le |\alpha| \le N = N(d).\}
\end{align}
for fixed constants $0 < M_1, M_2, \dots, M_N$. These estimates play a
key role in the large data theory for nontranslation-invariant
$L^2$-critical Cauchy problems typified
by~\eqref{e:intro_mass-crit_V}. We briefly discuss the nonlinear
problem in the last section of the introduction.


The case of one space dimension was treated in a previous joint work
with Killip and Visan~\cite{mkv_inv_str_1d}. This paper extends the
methods introduced there to higher dimensions.

\subsection{The setup}

To clarify the structure of our arguments we begin with a slightly
more general setup. Hence we consider time-dependent, real-valued
symbols $a(t, x, \xi)$ which are measurable in $t$ and satisfy
\begin{align}
  \label{e:intro_symbol_hyp}
  |\partial^\alpha_x \partial^\beta_\xi a| \le c_{\alpha \beta} \
  \text{for all} \ |\alpha| + |\beta| \ge 2.
\end{align}
Further, we assume the characteristic curvature condition
\begin{align}
  \label{e:char_curvature}
  \bigl| |\operatorname{det} a_{\xi \xi}|  - 1\bigr| +   \bigl| \| a_{\xi \xi} 
  \| - 1\bigr| \le \varepsilon
\end{align}
for some small $0 < \varepsilon < 1$.  For concreteness, all matrix norms in this article denote the
Hilbert-Schmidt norm, but the exact choice of norm is
inessential.

These hypotheses encompass several interesting situations:
\begin{itemize}
\item Schr\"{o}dinger Hamiltonians with time-dependent scalar
  potentials $a = \tfrac{1}{2}|\xi|^2 + V(t, x)$, where
  $V \in \mathcal{V}$.

\item Electromagnetic-type symbols
  $a = \tfrac{1}{2}|\xi|^2 + b(x, \xi) + V(t, x)$, where the first order
  symbol $b(x, \xi)$ is real and satisfies
  $|\partial^\alpha_x \partial^\beta_\xi b| \le c_{\alpha \beta}$ for
  all $|\alpha| + |\beta| \ge 1$, and $V \in \mathcal{V}$ is a scalar
  potential as before.
    
\item The frequency $1$ portion of the Laplacian on a curved
  background.
\end{itemize}

For a symbol as defined above, write $a^w(t, x, D)$ for its Weyl quantization. Let
$U(t, s)$ denote its unitary propagator on $L^2(\mathbf{R}^d)$, so that
$u := U(t, s) u_s$ is the solution to the equation
\begin{align}
  \label{e:intro_pd_evolution}
  (D_t + a^w(t, x, D)) u = 0, \quad u(s, \cdot) = u_s \in L^2(\mathbf{R}^d),
\end{align}
Evolution equations of this type were studied by Koch and
Tataru~\cite{KochTataru2005}. 
While translations and modulations do not preserve the
equation~\eqref{e:intro_pd_evolution}, they do preserve the class of
equations defined by our assumptions.  For an element $(x_0, \xi_0)$
of classical phase space, define the ``phase space translation''
operator $\pi(x_0, \xi_0)$ by
\begin{align}
\nonumber
 \pi(z_0)f (x) = e^{i \langle x-x_0, \xi_0 \rangle} f(x-x_0).
\end{align}
Then a direct computation, as in the proof
of~\cite[Proposition 4.3]{KochTataru2005}, yields
\begin{lemma}
  \label{l:galilei}
  If $U(t,s)$ is the propagator for the symbol $a$ and
  $\sigma \mapsto z^\sigma = (x^\sigma, \xi^\sigma)$ is a
  bicharacteristic of $a$, then
  \begin{align*}
    U(t, s) \pi(z_0^s) f = e^{i (\phi(t, z_0) - \phi(s, z_0))} \pi(z_0^t) 
U^{z_0}(t, s),
  \end{align*}
  where $U^{z_0}$ is the propagator for the equation
  \begin{gather*}
    [D_t + (a^{z_0})^w (t, x, D)] u = 0,\\
    a^{z_0} (t, z) = a(t, z_0^t + z) - \langle x,a_{x}(t,
   z_0^t)\rangle - \langle \xi ,a_{\xi}(t, z_0^t)\rangle - a(z_0^t),
  \end{gather*}
  and the phase is defined by
  \begin{align*}
    \phi(t, z_0) = \int_0^t \langle a_{\xi}(\tau, z_0^\tau),
    \xi_0^\tau\rangle - a(\tau, z_0^\tau) \, d\tau.
  \end{align*}
\end{lemma}
Observe that the transformed symbol $a^{z_0}$ satisfies the same
estimates assumed of $a$. As a special case, symbols of the form
$a = \tfrac{1}{2}|\xi|^2 + \langle A(t,x), \xi \rangle + \omega_{jk}(t)
x^j x^k$ are themselves preserved by the mapping $a \mapsto a^{z_0}$
if $A = A_j dx^j$ is a 1-form whose components are linear functions of
the space variables with time-dependent coefficients.  In two and
three space dimensions, such $A$ are potentials for uniform magnetic
fields.

The preceding hypotheses imply that the equation~\eqref{e:intro_pd_evolution} 
satisfies
a local-in-time dispersive estimate:
\begin{lemma}
  \label{l:intro_disp_est}
  If the symbol $a$ satisfies the
  conditions~\eqref{e:intro_symbol_hyp} and \ref{e:char_curvature}, there exists
  $T_0 > 0$ such that the propagator $U(t, s)$ for the evolution
  equation~\eqref{e:intro_pd_evolution} satisfies the estimate
  \begin{align}
    \label{e:intro_loc_disp}
    \| U(t,s) \|_{L^1_x \to L^\infty_x} \lesssim |t-s|^{-d/2} \
    \text{for all} \ |t-s| \le T_0.
  \end{align}
  Hence, the solutions to~\eqref{e:intro_pd_evolution} satisfy
  local-in-time Strichartz estimates
  \begin{align*}
    \|u\|_{L^q_t L^r_x( I \times \mathbf{R}^d)} \lesssim_{|I|} \| u_s 
    \|_{L^2(\mathbf{R}^d)}
  \end{align*}
  for any compact time interval $I$, and for all \emph{Strichartz
    exponents} $(q, r)$ satisfying $2 \le q, r \le \infty$,
  $\tfrac{2}{q} + \tfrac{d}{r} = \tfrac{d}{2}$, and
  $(q, r, d) \ne (2, \infty, 2)$.
\end{lemma}
\begin{proof}[Proof sketch]
  The dispersive estimate is shown in Koch-Tataru using wavepacket
  parametrices~\cite[Proposition 4.7]{KochTataru2005}. Standard
  arguments~(see Ginibre-Velo~\cite{Ginibre1995} and
  Keel-Tao~\cite{keel-tao}) then yield the Strichartz estimates.
\end{proof}
It suffices to choose the time increment $T_0$ so that
\begin{align}
  \label{e:time_smallness}
  T_0 \le 1, \quad T_{0} \| a_{x\xi}\| + T_0^2 \| a_{xx} \| \le \eta,
\end{align}
where $\eta = \eta(d)$ is a small parameter depending only on the
dimension.

\begin{remark}
  The concrete cases of scalar potentials and magnetic potentials were
  studied much earlier by Fujiwara and Yajima, respectively, who
  proved the dispersive bound using Fourier integral
  parametrices~\cite{fujiwara_fundamental_solution,Yajima1991}.
\end{remark}

We seek refinements of the Strichartz inequality analogous to those
for the constant-coefficient equation.  The earlier arguments for
constant coefficient equation relied crucially on subtle bilinear
estimates from Fourier restriction theory. We isolate and reformulate
the technical lynchpin in the present context.
\begin{hyp}
  \label{h:bilinear_Lp}
  There exist $T_0 > 0$ and $1 < p < \tfrac{d+2}{d}$ such
  that the following holds: if $f, g \in L^2(\mathbf{R}^d)$ have frequency
  supports in sets of diameter $\lesssim N$ which are separated by
  distance $\sim N$, then 
  \begin{align}
    \label{e:bilinear_Lp_hyp}
    \| U_{\lambda}^{s}(t) f U_{\lambda}^{s} (t) g \|_{L^p_{t,x} ([-T_0, T_0] \times
    \mathbf{R}^d)} \lesssim N^{-\delta} \|f\|_{L^2(\mathbf{R}^d)} \| 
    g\|_{L^2(\mathbf{R}^d)},
  \end{align}
  for all $s \in [-1, 1]$ and all $0 < \lambda \le 1$, where
  $U_{\lambda}^{s}(t) = U^s_\lambda (t, 0)$ are the propagators for
  the time-translated and rescaled symbols
  $a^{s}_{\lambda}:= \lambda^2 a(s + \lambda^2 t, \lambda x,
  \lambda^{-1} \xi)$.
\end{hyp}
When $a = \tfrac{1}{2}|\xi|^2$, the scaling and translation parameters
$\lambda, \ s$ are extraneous, and inequalities of the
form~\eqref{e:bilinear_Lp_hyp} are called (adjoint) bilinear Fourier
restriction estimates.  They were utilized by B\'{e}gout-Vargas to
obtain mass-critical Strichartz refinements in dimension $3$ and
higher~\cite{begout-vargas} (the results in dimensions 1 and 2, due to
Carles-Keraani, Merle-Vega, and Moyua-Vargas-Vega utilized linear
restriction estimates~\cite{carles-keraani,merle-vega,Moyua1999}). For
further discussion of such estimates see for instance \cite{Tao2003} and the
references therein.

In the first part of this paper, we connect~\eqref{e:bilinear_Lp_hyp} to
Strichartz refinements. To measure concentration in
the solution we test it against scaled, modulated, and translated
wavepackets.  Set
\begin{align}
  \label{e:gaussian_wp}
  \psi(x) = c_d e^{-\frac{|x|^2}{2}}, \ \psi_{x_0,\xi_0} = \pi(x_0,
  \xi_0) \psi, \quad c_d = 2^{-d/2} \pi^{-3d/4},
\end{align}
where $S_\lambda$ is the the unitary rescaling $S_\lambda f(x) := 
\lambda^{-d/2} f(\lambda^{-1} x)$.
\begin{theorem}
  \label{t:conditional_inv_str}
  If Hypothesis~\ref{h:bilinear_Lp} holds, then there exists
  $0 < \theta < 1$ such that for all initial data
  $u_0 \in L^2(\mathbf{R}^d)$, the solution $u$ to the
  equation~\eqref{e:intro_pd_evolution} satisfies
\begin{align}
  \label{e:refined_str}
  \|u\|_{L^{\frac{2(d+2)}{d}} ([-1,1] \times \mathbf{R}^d)} \lesssim 
  \bigl(\sup_{0
  < \lambda \le 1, \ |t| \le 1, \ (x_0, \xi_0) \in T^* \mathbf{R}^d} |\langle 
  S_\lambda
  \psi_{x_0, \xi_0}, u(t) \rangle_{L^2(\mathbf{R}^d)}| \bigr)^{\theta} \|
  u_0\|_{L^2(\mathbf{R}^d)}^{1-\theta}.
\end{align}

\end{theorem}

The generality of our hypotheses requires us to formulate the
estimates locally in time. Indeed, for most potentials the left side
of the Strichartz estimate~\eqref{e:refined_str} is infinite if one
integrates over $\mathbf{R} \times \mathbf{R}^d$; for instance, the harmonic
oscillator potential $V=|x|^2$ admits periodic-in-time
solutions. Nonetheless, our methods do yield (a new proof of) a
global-in-time refined Strichartz estimate
\begin{align*}
  \|u \|_{L^{\frac{2(d+2)}{d}}_{t,x}(\mathbf{R} \times
  \mathbf{R}^d)} \lesssim \Bigl( \sup_{\lambda > 0, \ t \in \mathbf{R}, \
  (x_0, \xi_0) \in T^* \mathbf{R}^d} | \langle S_\lambda \psi_{x_0,
  \xi_0}, u(t) \rangle_{L^2(\mathbf{R}^d)} | \Bigr)^{\theta} \| u_0
  \|_{L^2(\mathbf{R}^d)}^{1-\theta}.
\end{align*}
for solutions to the constant coefficient
equation~\eqref{e:intro_free_particle}.

In applications to PDE, such a refined estimate is nowadays interpreted
in the framework of concentration compactness, and yields profile
decompositions via repeated application of the following
\begin{lemma}
  \label{l:profile}
  Assume the estimate~\eqref{e:refined_str} holds. Let $u_n := U(t) f_n$
  be a sequence of linear solutions with initial data
  $u_n(0) = f_n \in L^2(\mathbf{R}^d)$ such that
  $\|f_n\|_{L^2(\mathbf{R}^d)} \le A < \infty$ and
  $\|u_n\|_{L^{\frac{2(d+2)}{d}}_{t,x}} \ge \varepsilon > 0$. Then,
  after passing to a subsequence, there exist parameters
  \begin{align*}
    \{ (\lambda_n, t_n, x_n, \xi_n) \}_n \subset (0, 1] \times [-1, 1]
    \times \mathbf{R}^d_x \times \mathbf{R}^d_\xi
  \end{align*}
  and a function $0 \ne \phi \in L^2(\mathbf{R}^d)$ such that
  \begin{gather*}
    \pi(x_n, \xi_n)^{-1} S_{\lambda_n}^{-1} u_n \rightharpoonup \phi \
    \text{in} \ L^2\\
    \| \phi\|_{L^2} \gtrsim \varepsilon \Bigl(\frac{\varepsilon}{A}
    \Bigr)^{\frac{1-\theta}{\theta}}.
  \end{gather*}
  Further,
  \begin{align*}
    \|f_n \|_{L^2}^2 - \|f_n - U(t_n)^{-1} S_{\lambda_n} \pi(x_n, \xi_n) S_{\lambda_n}
    \phi\|_{L^2}^2 - \| U(t_n)^{-1} S_{\lambda_n} \pi(x_n, \xi_n) S_{\lambda_n}
    \phi\|_{L^2}^2 \to 0.
  \end{align*}
\end{lemma}

\begin{proof}
  By the estimate~\eqref{e:refined_str}, there exist $\lambda_n, t_n,
  x_n, \xi_n$ such that
  \begin{align*}
    |\langle S_{\lambda_n} \psi_{x_n, \xi_n}, U(t_n)f_n \rangle|  =
    |\langle \psi, \pi(x_n, \xi_n)^{-1} S_{\lambda_n}^{-1} U(t_n) f_n \rangle| \gtrsim
    \varepsilon \Bigl( \frac{\varepsilon}{A} \Bigr)^{\frac{1-\theta}{\theta}}.
  \end{align*}
  The sequence $\pi(x_n, \xi_n)^{-1} S_{\lambda_n}^{-1} U(t_n) f_n$ is
  bounded in $L^2$, and therefore converges weakly in $L^2$ to some
  $\phi$ after passing to a subsequence. The lower bound on $\|\phi\|_{L^2}$ is
  immediate, while
  \begin{align*}
    &\|f_n\|_{L^2}^2 - \|f_n - U(t_n)^{-1} S_{\lambda_n} \pi(x_n, \xi_n)
    \phi\|_{L^2}^2 - \| U(t_n)^{-1}S_{\lambda_n} \pi(x_n, \xi_n)
       \phi\|_{L^2}^2\\
    &= 2 \operatorname{Re} \langle f_n - U(t_n)^{-1}S_{\lambda_n}  \pi(x_n, 
    \xi_n)
     \phi, U(t_n)^{-1} S_{\lambda_n} \pi(x_n, \xi_n) 
      \phi \rangle\\
    &= 2 \operatorname{Re} \langle  \pi(x_n, \xi_n)^{-1} S_{\lambda_n}^{-1}
      U(t_n) f_n - \phi, \phi \rangle \to 0.
  \end{align*}
\end{proof}
Further discussion of profile decompositions and inverse Strichartz theorems 
may be found in the lecture notes~\cite{claynotes} and the references therein.

In the second part of this paper, we verify
Hypothesis~\ref{h:bilinear_Lp} for scalar potentials.

\begin{theorem}
  \label{t:bilinear_Lp}
  Consider a Schr\"{o}dinger operator of the form
  $H(t) = -\tfrac{1}{2}\Delta + V(t,x)$, where $V \in \mathcal{V}$. Suppose
  $S_1, S_2 \subset \mathbf{R}^d_{\xi}$ are subsets of Fourier space with
  $\operatorname{diam}(S_j) \le N$ and $c^{-1}N \ge \operatorname{dist}(S_1, 
  S_2) \ge cN$ for some
  $0 < c < 1$.  There exists a constant $\eta = \eta(c) > 0$ such that
  if $\tau_0 > 0$ satisfies
  \begin{align*}
    (\tau_0 + \tau_0^2) \| \partial^2_xV\|_{L^\infty} < \eta,
  \end{align*}
  then for any $f, g \in L^2(\mathbf{R}^d)$ with
  $\operatorname{supp}(\hat{f}) \subset S_1$ and
  $\operatorname{supp}(\hat{g}) \subset S_2$, the corresponding
  linear solutions $u= U(t, 0) f$ and $v = U(t, 0) g$
  satisfy the estimate
  \begin{align}
    \label{e:bilinear_Lp}
    \| u v\|_{L^q ( [-T_0 , T_0 ] \times \mathbf{R}^d)}
    \lesssim_{\varepsilon} N^{d - \frac{d+2}{q} +\varepsilon} \|f\|_{L^2} 
    \|g\|_{L^2} \quad
    \text{for all} \quad \frac{d+3}{d+1} \le q < \frac{d+2}{d}
  \end{align}
  for any $\varepsilon > 0$, $N \ge 1$, and $V \in \mathcal{V}$.
\end{theorem}
For $V=0$, the above estimate was conjectured by Klainerman and
Machedon without the epsilon loss, and first proved by Wolff for the
wave equation~\cite{Wolff2001} and subsequently by Tao~\cite{Tao2003}
for the Schr\"{o}dinger equation (both with the epsilon
loss). Strictly speaking, the time truncation is not present in the
original formulations of those estimates, but may be easily removed by
a rescaling and limiting argument.

Finally, while we make no attempt to address general
magnetic potentials, a simple case with some physical relevance
does essentially follow from the proof for scalar potentials. The
necessary modifications for the following theorem are sketched in the last section.
\begin{theorem}
  \label{t:bilinear_Lp_magnetic}
  The conclusion of the previous theorem holds for Schr\"{o}dinger 
operators of the form $H(t) =
  -\tfrac{1}{2}(\nabla - iA )^2 + V(t, x)$ where $A = A_j dx^j$ is a
  1-form whose components are linear in the space variables (i.e. the
  vector potential for a uniform magnetic field), and condition on the
  time increment $\tau_0$ is replaced by
  \begin{align*}
    \tau_0 \| a_{x \xi}\| + (\tau_0 + \tau_0^2)\| a_{xx}\| < \eta.
  \end{align*}
\end{theorem}

We remark that the restriction estimate~\eqref{e:bilinear_Lp_hyp} does
\emph{not} hold for all symbols satisfying the
conditions~\eqref{e:intro_symbol_hyp}
and~\eqref{e:char_curvature}. For instance, it was observed by
Vargas~\cite{Vargas2005} that when
$U(t) = e^{it \partial_x \partial_y}$ is the ``nonelliptic''
Schr\"{o}dinger propagator in two space dimensions (thus
$a = \xi_x \xi_y$), the bilinear restriction
estimate~\eqref{e:intro_symbol_hyp} can fail unless the frequency
supports of the two inputs are not only disjoint but also separated in
both Fourier coordinates. In fact, the
refinement~\eqref{e:refined_str} as stated is \emph{false} for the
nonelliptic equation; for a correct formulation, one should enlarge
the symmetry group on the right side to include the hyperbolic
rescalings $u(x, y) \mapsto u(\mu x, \mu^{-1} y)$; see the work of
Rogers and Vargas~\cite{Rogers2006a}.

While the classical bicharacteristics of elliptic and nonelliptic
propagators seemingly have no qualitative difference---and indeed the
dispersive estimates hold equally well for both---the quantum
propagators have radically different behavior in terms of oscillations
in time. If one compares the travelling wave solutions
\[
  e^{i [x \xi_x + y \xi_y - \frac{t}{2}(\xi_x^2 +
    \xi_y^2)]}, \quad   e^{i [ x \xi_x + y \xi_y -
    t \xi_x \xi_y]},
\]
it is evident that unlike in the elliptic case, two solutions to the
nonelliptic equation which are well-separated in spatial frequency
need not decouple in time.

The lesson of this counterexample is that while the dispersive and
Strichartz estimates follow directly from properties of the classical
Hamiltonian flow, an \emph{inverse} Strichartz estimate depends more
subtly on the temporal oscillations of the quantum evolution, which is
connected to the bilinear decoupling estimates.

\subsection{The main ideas}

Suppose one has initial data $u_0 \in L^2$ such that the corresponding
solution $u$ has nontrivial Strichartz norm. Then, we need to identify
a bubble of concentration in $u$, characterized by several parameters
that reflect the underlying symmetries in the problem. In the
$L^2$-critical setting, the relevant features consist of a
significant length scale $\lambda_0$ as well as the position $x_0$,
frequency $\xi_0$, and time $t_0$ when concentration occurs.

The existing proofs of Strichartz refinements for the
constant-coefficient equation first use spacetime Fourier analysis
(including restriction estimates) to identify a cube $Q$ in Fourier
space accounting for a significant portion of the spacetime norm of
$u$, which reveals the frequency center $\xi_0$ and scale $\lambda_0$
of the concentration. For example, Begout-Vargas~\cite{begout-vargas}
first establish an extimate of the form
\begin{align*}
  \|e^{\frac{it\Delta}{2}} f\|_{L^{\frac{2(d+2)}{d}}} \lesssim \Bigl( \sup_{Q \
  \text{dyadic cubes}}
  |Q|^{1 - \frac{p}{2}} \int_{Q} |\hat{f}(\xi)|^p \, d\xi \Bigr)^{\mu}
  \|f\|_{L^2(\mathbf{R}^d)}^{1 - \mu p}
\end{align*}
Then, the time $t_0$ and position $x_0$ are recovered via a separate
physical-space argument. These arguments ultimately rely on
the fact that when $V = 0$, the equation is diagonalized by the
Fourier transform.

For equations with variable coefficients, it is more natural to
consider position $x_0$ and frequency $\xi_0$ together as a point in
phase space, which propagates along the bicharacteristics for the
equation. Following the approach in~\cite{mkv_inv_str_1d} for the
one-dimensional equation, we work in the physical space and first
isolate a significant time interval
$[t_0 - \lambda_0^2, t_0 + \lambda_0^2]$, which also suggests a
characteristic scale $\lambda_0$. Then $x_0$ and $\xi_0$ are recovered
by phase space techniques.

The first part of the argument in~\cite{mkv_inv_str_1d} carries over
essentially unchanged; however, the ensuing phase space analysis in
higher dimensions is more involved and occupies the bulk of this
article.

\subsection{An application to mass-critical NLS}

This article was originally motivated by the problem of proving global
wellposedness for the mass-critical quantum harmonic oscillator
\begin{align}
  \label{e:mcrit_harmonic}
  i\partial_t u = \Bigl( -\frac{1}{2}\Delta + \sum_{j} \omega_j^2 x_j^2
  \Bigr) u \pm |u|^{\frac{4}{d}} u.
\end{align}

By spectral theory, the Cauchy problem for~\eqref{e:mcrit_harmonic} is naturally posed in
the ``harmonic'' Sobolev spaces
\[u_0 \in \mathcal{H}^s := \{ u_0 \in L^2 : (-\Delta + \sum_{j}
  \omega_j^2 |x|^2)^{s/2}, \ u_0 \in L^2\}\] Global existence for data
in the ``energy'' space $\mathcal{H}^1$ was studied by
Zhang~\cite{Zhang2005}. More recently, Poiret, Robert, and Thomann
established probabilistic wellposedness in two space dimensions for
all subcritical cases $0 < s < 1$, as well as for other supercritical
problems~\cite{Poiret2014}. Another recent contribution by Burq,
Thomann, and Tzvetkov constructs Gibbs measures and proves
probabilistic global wellposedness for the critical case in one
dimension~\cite{Burq2013}.

It is well-known that the \emph{isotropic} harmonic oscillator
$\omega_j \equiv \tfrac{1}{2}$ may be ``trivially'' solved; to construct
solutions on unit length time intervals for arbitrary $L^2$ data, it
suffices to observe that $u$ is a solution
of~\eqref{e:intro_mass-crit_nls} on $\R_t \times \R_x^d$ iff its
\emph{Lens transform}
\begin{align*}
  \mathcal{L} u (t, x) : = \frac{1}{(\cos t)^{d/2}} u \Bigl( \tan t,
  \frac{x}{\cos t} \Bigr) e^{-\frac{i|x|^2 \tan t}{2}}
\end{align*}
solves~\eqref{e:mcrit_harmonic} on $(-\pi/2 \times 
\pi/2)_t \times \R_x^d$ with the same initial data. However, this trick
relies on algebraic cancellations that no longer hold for more general
harmonic oscillators. For further discussion of the nonlinear harmonic oscillator
as well as its connection with the Lens transform, consult the article of
Carles~\cite{carles_time-dependent_potential}.

To solve~\eqref{e:mcrit_harmonic} for large data in the critical space
$L^2$, the concentration compactness and rigidity approach is much
more promising.  Experience has shown that constructing suitable
profile decompositions is a core difficulty implementing this strategy
for dispersive equations with broken symmetries (e.g. loss of
translation-invariance). For instance,
see~\cite{me_quadratic_potential} for the energy-critical variant of
the quantum harmonic oscillator, as well
as~\cite{IonescuPausaderStaffilani2012,kvz_exterior_convex_obstacle},
and the references therein, for other energy-critical NLS on
non-Euclidean domains. Thus this article supplies the main harmonic
analysis input for the deterministic large data theory
of~\eqref{e:mcrit_harmonic} at the critical regularity.

\textbf{Acknowledgements}. The author is grateful to Michael Christ,
Rowan Killip, Daniel Tataru, and Monica Visan for many helpful
discussions, and also wishes to thank the anonymous referee for
numerous suggestions for improving the original manuscript. This
research was partially supported by the National Science Foundation
under Award No. 1604623. Part of this work was completed during the
2017 Oberwolfach workshop in ``Nonlinear Waves and Dispersive Equations.''

\section{Preliminaries}

\subsection{Notation}
We use the Japanese bracket notation $\langle x \rangle :=
(1+|x|^2)^{\frac{1}{2}}$.

\subsection{Classical flow estimates}

 We collect some elementary properties of the classical Hamiltonian
flow
\begin{align}
  \label{e:hamilton_ode}
  \left\{\begin{array}{ll} \dot{x} = a_{\xi}, & x(0) = y\\
           \dot{\xi} = -a_x, & \xi(0) = \eta.
                               \end{array}\right.
\end{align}
Solutions to this system are \emph{bicharacteristics}. For a point
$z = (x, \xi)$ in phase space, let
$\sigma \mapsto z^\sigma = (x^\sigma, \xi^\sigma)$ denote the
bicharacteristic initialized at $(x, \xi)$. Write
$(y, \eta) \mapsto (x^t(y, \eta), \xi^t(y, \eta))$ for the flow map.

The linearization of~\eqref{e:hamilton_ode} satisfies the following
Gronwall estimates:
\begin{lemma}
  \label{l:bichar}
  Suppose $|t| \|\partial_{x,\xi}^2 a\|_{L^\infty} \le 1$. Then
  \begin{align}
    \begin{split}
    \frac{\partial x^t}{\partial \eta} &= \int_0^t a_{\xi \xi}(\tau, x^\tau, 
\xi^\tau) \, d\tau +
    O(t^2 \|a_{x\xi} \| a_{\xi \xi}\|) + O( t^3 \|a_{xx} \| \| a_{\xi
      \xi} \|^2)\\
    \frac{\partial \xi^t}{ \partial \eta} &= I + O( t\|a_{\xi x} \|) +
                                          O( t^2 \|a_{xx}\| \|
                                          a_{\xi\xi} \|)\\
    \frac{\partial x^t}{\partial y} &= I + O(t\|a_{x\xi}\|) + O(
                                           t^2 \| a_{xx} \| \| a_{\xi
                                           \xi} \|)\\
    \frac{\partial \xi^t}{\partial y} &= \int_0^t- a_{xx}(\tau, x^\tau, 
    \xi^\tau) 
\, d\tau + O(
                                      t^2\|a_{xx} \| \| a_{x \xi} \|) +
                                      O( t^3 \|a_{xx} \|^2 \|a_{\xi
                                        \xi}\|)
                                    \end{split}.
  \end{align}
\end{lemma}

\begin{proof}
  The linearized system takes the form
  \begin{align*}
    \dot{y} &= a_{\xi x} y + a_{\xi \xi} \eta, \\
    \dot{\eta} &= -a_{xx}y - a_{x\xi} \eta.
  \end{align*}
  A preliminary application of Gronwall implies
  $|y(t)| + |\eta(t)| \lesssim |y(0)| + |\eta(0)|$.

  Consider  initial data $y(0) = I, \ \eta(0) = 0$.
  Then
  \begin{align*}
    |\eta(t)| \le \int_0^t |a_{xx} y| \, d\tau + \int_0^t |a_{x\xi}
    \eta(\tau)| \, d\tau,
  \end{align*}
  so
  $
    |\eta(t)| \lesssim  t\| a_{xx} \|.
  $
  Substituting this into the equation for $y$, we deduce
  \begin{align*}
    |y - I| \le \int_0^t |a_{\xi x} y| \, d\tau + \int_0^t
    |a_{\xi\xi} \eta| \, d\tau \lesssim t \| a_{\xi x} \| + t^2 \|
    a_{\xi \xi} \| \| a_{xx} \|.
  \end{align*}
  This in turn yields the refinement
  \begin{align*}
    \Bigl |\eta(t)  + \int_0^t a_{xx} \, d\tau \Bigr| \lesssim t^2 \|
    a_{xx} \| \| a_{\xi x} \| + t^3 \| a_{xx} \|^2 \| a_{\xi \xi} \|.
  \end{align*}

  The case $y(0) = 0, \ \eta(0) = I$ is similar. We have
  \begin{align*}
    |y(t)| \le \int_0^t |a_{\xi \xi} \eta| \, d\tau + \int_0^t |a_{\xi
    x} y| \, d\tau \Rightarrow |y(t)| \lesssim t \|a_{\xi \xi}\|,
  \end{align*}
  which yields
  \begin{align*}
    |\eta(t) - I| &\lesssim  \int_0^t \| a_{xx} \| \| a_{\xi \xi} \|
    \tau \, d\tau + \int_0^t | a_{x\xi} \eta| \, d\tau \lesssim t \|
                    a_{x\xi} \| + t^2 \| a_{xx} \|  \| a_{\xi\xi}
                    \|,\\
    \Bigl| y(t) - \int_0^t a_{\xi \xi } \, d\tau \Bigr| &\lesssim  t^2
    \| a_{\xi x} \| \| a_{\xi\xi} \| + t^3 \| a_{xx} \| \| a_{\xi \xi}\|^2.
  \end{align*}
\end{proof}


These imply, in view of the normalizations~\eqref{e:char_curvature},
the integrated estimates
\begin{gather}
  \label{e:integrated_bichar}
  \begin{split}
    &\begin{split}
      x^t_1 - x_2^t &= x_1^s - x_2^s + [I' + O(\varepsilon)] (t-s) (\xi_1^s -
      \xi_2^s)\\
      &+ O(|t-s| \|a_{x\xi}\|) (|x_1^s - x_2^s| + |t-s|
      |\xi_1^s - \xi_2^s|)\\
      &+ O(|t-s|^2 \| a_{xx}\| )(|x_1^s - x_2^s|) +
      |t-s|   |\xi_1^s - \xi_2^s|).
    \end{split}\\
    &\begin{split}
      \xi_1^t - \xi_2^t &= \xi_1^s - \xi_2^s \\
      &+ O(
      |t-s|\|a_{xx}\|) |x_1^s - x_2^s| \\
      &+ O(|t-s|^2 \|a_{xx} \| \|a_{x\xi}\|) |x_1^s - x_2^s|
      +O(|t-s| \| a_{x \xi}\|) |\xi_1^s - \xi_2^s|\\
      &+ O(|t-s|^3 \| a_{xx} \|^2) |x_1^s - x_2^s| + O(|t-s|^2
      \|a_{xx}\|) |\xi_1^s - \xi_2^s|,
    \end{split}
  \end{split}
\end{gather}
where $I'$ is an orthogonal matrix which equals the identity if
$a_{\xi\xi}$ is positive-definite.  In particular, we have
\begin{corollary}
  \label{c:once_collision}
  If $|x_1^s - x_2^s| \le r$, then $|x_1^t - x_2^t| \ge Cr$ whenever
  $\tfrac{ 2Cr}{|\xi_1^s - \xi_2^s|} \le |t-s| \le T_0$.
\end{corollary}
Physically, this means that two particles colliding with sufficiently
large relative velocity will only interact once in the time window of
interest.


Next, we record a technical lemma first proved in the 1d
case~\cite[Lemma 2.2]{mkv_inv_str_1d}. This is used in the proof of
Lemma~\ref{l:preimage_measure} below but the computations use the
preceding estimates.
\begin{lemma}
	\label{l:technical_lma}
	There exists a constant $C = C(\|\partial^2a\|) > 0$ so that if
	$Q_\eta = (0, \eta) +[-1, 1]^{2d} \subset T^* \mathbf{R}^d$ and $r \ge 1$, 
	then
	\[
	\bigcup_{|t - t_0| \le \min(|\eta|^{-1}, 1) } \Phi(t)^{-1} (z_0^t + 
	rQ_\eta) \subset
	\Phi(t_0)^{-1}  (z_0^{t_0}+ Cr Q_\eta).
	\]
\end{lemma}
In other words, if the bicharacteristic $z^t$ starting at $z \in 
T^*\mathbf{R}^d$ passes
through the cube $z_0^t + rQ_{\eta}$ in phase space during some time
window $|t - t_0| \le \min(|\eta|^{-1},1)$, then it must lie in the
dilate $z_0^{t_0} + CrQ_{\eta}$ at time $t_0$.

\begin{proof}
  If $z \in \Phi(t)^{-1} ( z_0^t + rQ_\eta)$, by definition we have
  $|x^t - x_0^t| \le r$ and $|\xi^t - \xi_0^t - \eta| \le r$. Assuming
  that $|\eta| \ge 1$, the estimates~\eqref{e:integrated_bichar} imply
  that
  \begin{align*}
    |x^{t_0} - x_0^{t_0}| &\le r + |\eta|^{-1}(|\eta| + r) \\
                          &+ O(|\eta|^{-1} \| \partial^2a\|)
                            ( r + |\eta|^{-1}( |\eta| + r)) + O(
                            |\eta|^{-2} \| \partial^2 a\|) ( r +
                            |\eta|^{-1} ( |\eta| + r))\\
                          &\le Cr\\
    |\xi^{t_0} - \xi_0^{t_0} - \eta| &\le r + O(|\eta|^{-1} \| a_{xx}
    \|) r + (|\eta|^{-2} \| a_{xx} \| \| a_{x\xi} \|) r + O(
    |\eta|^{-1} \| a_{x\xi} \|) (|\eta| + r)\\
    &+ (|\eta|^{-3} \| a_{xx}\|^2) r + O( |\eta|^{-2} \|a_{xx}|) (
      |\eta| + r)\\
    &\le Cr.
  \end{align*}
  The case $|\eta| < 1$ is similar.
\end{proof}

\subsection{Wavepackets}
\label{ss:wp_decomp}

Let $R \ge 1$ be a scale and $z_0 = (x_0, \xi_0)$ be a point in phase space.  A scale-$R$ \emph{wavepacket} at $z_0$ is a
Schwartz function $\phi_{z_0}$ such that $\phi_{z_0}$ and its
Fourier transform $\widehat{\phi_{z_0}}$ concentrate in the regions
$|x-x_0| \le R^{1/2}$ and $|\xi - \xi_0| \le R^{-1/2}$, respectively:
\begin{align*}
  |(R^{1/2} \partial_x)^k \phi_{z_0}(x)| \lesssim_{k,N} \Bigl \langle
  \frac{x-x_0}{R^{1/2}} \Bigr \rangle^{-N}, \quad   |(R^{-1/2}
  \partial_\xi)^k \widehat{\phi_{z_0}}(\xi)| \lesssim_{k,N} \Bigl \langle
  \frac{\xi-\xi_0}{R^{-1/2}} \Bigr\rangle^{-N} \quad \forall k, N\ge 0.
\end{align*}
There are many ways to decompose $L^2$ functions into linear
combinations of wavepackets. For the first part of this article, it is
technically more convenient to use a continuous decomposition. Later
on in Section~\ref{ss:discrete_wp_decomp}, we switch to a discrete
version which is more common in the restriction theory literature.

In this section we recall a standard continuous wavepacket transform.  To
keep things simple we work at unit scale since that is all we shall
need. For a function $f \in L^2(\mathbf{R}^d)$, its \emph{Bargmann
  transform} or \emph{FBI transform} is the function
$Tf \in L^2(T^* \mathbf{R}^d)$ defined by
\begin{align*} Tf(z) = \langle
  f, \psi_{z} \rangle_{L^2(\mathbf{R}^d)}, \quad \psi_z = \pi(z) \psi \text{
  as in } \eqref{e:gaussian_wp}.
\end{align*}
The transform satisfies a Plancherel identity
$\| Tf\|_{L^2( T^* \mathbf{R}^d)} = \|f\|_{L^2(\mathbf{R}^d)}$; dually,
for any wavepacket coefficients $F  \in L^2(T^*\mathbf{R}^d_z)$,
one has
\begin{align*}
  \|T^* F\|_{L^2_x} = \Big\| \int_{T^* \mathbf{R}^d} F(z) \psi_z
  \, dz \Bigr\|_{L^2_x} \le \|F\|_{L^2_z}.
\end{align*}
Indeed, $TT^*$ is the orthogonal projection onto $TL^2(\mathbf{R}^d)$. Then as 
$T^*T = I$, any $f \in L^2(\R^d)$ can be resolved (nonuniquely) into a 
continuous superposition of wavepackets \[f(x) = \int_{T^*\R^d} f_z \psi_z 
(x)\, dz.\]

Applying the propagator $U(t)$ to both sides and using 
linearity and the next lemma,  one obtains a wavepacket
decomposition
\[ u(t, x) = \int u_z(t, x) \, dz, \ u_z(t, x) = f_z [U(t) \psi_z] (x)\]
of
Schr\"{o}dinger solutions. For brevity we sometimes omit the arguments
and write $f = \int f_z \, dz, \ u = \int u_z \, dz$.

\begin{lemma}[Evolution of a packet]
  \label{l:wp_propagation}
  If $\psi_{z_0}$ is a scale-$1$ wavepacket, $U(t)$ is the propagator
  for the equation~\eqref{e:intro_pd_evolution}, and $z_0 \mapsto z_0^t$ is
  the bicharacteristic starting at $z_0$, then $U(t) \psi_{z_0}$ is a
  scale-$1$ wavepacket concentrated at $z_0^t$ for all
  $|t| = O(1)$.
\end{lemma}
\begin{proof}[Proof sketch]
  Using Lemma~\ref{l:galilei} we reduce to the case $z_0 = 0$ and also
  ensure that the symbol $a(t,x, \xi)$ vanishes to second order at
  $(x, \xi) = (0, 0)$ in addition to satisfying the
  bounds~\eqref{e:intro_symbol_hyp}. Then it suffices to show that
  propagator $U(t)$ for such symbols maps Schwartz functions to
  Schwartz functions on unit time scales. This is done using weighted
  Sobolev estimates as in~\cite[Section~4]{KochTataru2005}.
\end{proof}
The term \emph{wavepacket} shall also refer to spacetime functions of
the form $U(t) \psi_{z}$, not just the fixed time slices. Later it
will be essential to exploit not just the spacetime localization of
wavepackets but also their phase as described in
Lemma~\ref{l:galilei}.

\section{Choosing a length scale}
\label{s:time_scale}

We begin with the following lemma from~\cite[Proposition
3.1]{mkv_inv_str_1d}, obtained by a variant of the usual $TT^*$
derivation of the Strichartz estimates. While that article concerned
just Schr\"{o}dinger operators with scalar potentials, the proof works
equally well in the current more general setting.
\begin{proposition}
  Suppose $U(t,s)$ satisfies a local in time dispersive estimate as in
  Lemma~\ref{l:intro_disp_est}. Let $(q, r)$ be Strichartz exponents
  (i.e. satisfying the conditions in that Lemma) with
  $2 < q < \infty$. Assume that  $f \in L^2(\mathbf{R}^d)$ satisfies
  $\|f\|_{L^2(\mathbf{R}^d)} = 1$ and
  \begin{align*}
    \| U(t) f\|_{L^{q}_t L^r_x ( [-1, 1] \times \mathbf{R}^d )}
    \ge \varepsilon.
  \end{align*}
  Then there is a time interval $J \subset [-1, 1]$ such that
  \begin{align*}
    \| U(t, s)f\|_{L^{q-1}_t L^r_x( J \times \mathbf{R}^d)} \gtrsim
    |J|^{\frac{1}{q(q-1)}} \varepsilon^{ \frac{q}{q-2} }.
  \end{align*}
  Equivalently, 
  \begin{align*}
    \| U(t, s) f\|_{L^{q} L^r} \lesssim \Bigl( \sup_{J \in [-1, 1]} 
    |J|^{-\frac{1}{q(q-1)}} \| U(t,s)f\|_{L^{q-1}_t L^r_x ( J \times
    \mathbf{R}^d)} \Bigr)^{1 - \frac{2}{q}} 
    \|f\|_{L^2(\mathbf{R}^d)}^{\frac{2}{q}}.
  \end{align*}
\end{proposition}
Note that by pigeonholing we may always assume that $|J| \le T_0$,
where $T_0$ is the time increment selected
in~\eqref{e:time_smallness}.

Now let $(q, r)$ be the Strichartz exponents determined by the
conditions $\tfrac{2}{q} + \tfrac{d}{r} = \tfrac{d}{2}$ and $q-1 = r$. It is
easy to see that $2 < r < \tfrac{2(d+2)}{d} < q < \infty$.

For each $J = [s - \mu, s+ \mu] \subset [-1, 1]$, we write
\[
U(t, s) f =   \Bigl(
\frac{T_0}{\mu} \Bigr)^{d/4} \tilde{U} \Bigl( \frac{T_0}{\mu} (t-s),
0\Bigr)\tilde{f} \Bigl ( \sqrt{ \frac{T_0}{\mu} } x), \quad \tilde{f} =  \Bigl( 
\frac{\mu}{T_0}
  \Bigr)^{d/4} f\Bigl( \sqrt{\frac{\mu}{T_0}} x \Bigr),
\]
where $\tilde{U}(t, s)$ is the propagator for the rescaled equation
$(D_t + \tilde{a}^w ) \tilde{u} = 0$, and
\[\tilde{a}(t, x, \xi) := \frac{\mu}{T_0} a\Bigl( s + \frac{\mu}{T_0} t,
\sqrt{\frac{\mu}{T_0}} x, \sqrt{\frac{T_0}{\mu}} \xi \Bigr).\]
Changing variables, we obtain
\begin{align*}
  |J|^{-\frac{1}{q(q-1)}} \| U(t, s) f\|_{L^{q-1}_t L^r (J \times
  \mathbf{R}^d)} = \| \tilde{U}(t) \tilde{f} \|_{L^{q-1}_t L^r_x([-T_0,
  T_0] \times \mathbf{R}^d)}.
\end{align*}

By interpolating with $L^2_{t,x} ([-T_0, T_0] \times \mathbf{R}^d)$, which
is bounded by unitarity, we see that
Theorem~\ref{t:conditional_inv_str} would follow if we prove that for
some $2 < q_0 < \tfrac{2(d+2)}{d}$ and $0 < \theta < 1$, the scale-1
refined estimate
\begin{align}
  \label{e:refined_Lp}
  \| U_\lambda^s(t) f\|_{L^{q_0}([-T_0, T_0] \times \mathbf{R}^d)} \lesssim
  (\sup_{z} | \langle \psi_{z}, f \rangle|)^\theta  \|f\|_{L^2}^{1-\theta}.
\end{align}
holds for all $s \in [-1, 1], \ 0 < \lambda \le 1$, where the notation
$U_{\lambda}^s(t)$ is as in Hypothesis~\ref{h:bilinear_Lp}.

Over the next two sections we establish
\begin{proposition}
  \label{p:reduced_thm1}
  If Hypothesis~\ref{h:bilinear_Lp} holds, then so does the
  estimate~\eqref{e:refined_Lp}.
\end{proposition}

\section{A refined bilinear $L^2$ estimate}
\label{s:refined_L2}

In previous work~\cite{mkv_inv_str_1d}, we proved~\eqref{e:refined_Lp}
when $d=1$ with $q_0 = 4$ by viewing the inequality as a bilinear
$L^2$ estimate and exploit orthogonality. Such a direct approach fails
in $d\ge 2$ dimensions; since $2 < \frac{2(d+2)}{d} \le 4$, the
left side of~\eqref{e:refined_Lp} could well be infinite when
$q_0 = 4$. To obtain a refined linear $L^{q_0}$ estimate for
$q_0 < \tfrac{2(d+2)}{d}$, we also begin by interpreting it as a refined
bilinear $L^{q_0/2}$ estimate, but use dyadic decomposition and
interpolate between two microlocalized estimates:
\begin{itemize}
\item A refined bilinear $L^2$ estimate (``refined'' in the
  sense of exhibiting a sup over wavepacket coefficients) with some loss in
  the frequency separation of the inputs.
\item A bilinear $L^p$ estimate  for some $p < \tfrac{d+2}{d}$ which
  yields gains in the frequency separation, essentially the content of 
Hypothesis~\ref{h:bilinear_Lp}.
\end{itemize}

This section discusses the former. In the next section we put together
the two estimates, and the $L^p$ estimate is established in the
remainder of the paper.

\begin{proposition}
  \label{p:dy_refined_L2}
  Suppose $f = \int f_z \psi_z \, dz$ and $g = \int g_z \psi_z \, dz$ are
  $L^2(\mathbf{R}^d)$ initial data with corresponding Schr\"{o}dinger
  evolutions $u = \int u_z \, dz$ and $v = \int v_z \, dz$, where $u_z(t,x) = 
f_z[ U(t)
\psi_z](x)$, $v_z (t,x)= g_z [U(t) \psi_z] (x)$. Then 
  \begin{align}
    \label{e:dy_refined_L2}
    \Bigl\|\int_{|\xi_1 - \xi_2| \sim N} u_{z_1} v_{z_2} \,
    dz_1 dz_2 \Bigr \|_{L^2([-T_0, T_0]\times \mathbf{R}^d)} \lesssim N^\alpha
    (\sup_{z} |f_{z}|^{1/p'} \| f_{z} \|_{L^2_z}^{1/p})(
    \sup_{z} |g_{z}|^{1/p'} \| g_{z} \|_{L^2_{z}}^{1/p})
  \end{align}
  for some $\alpha = \alpha(d) $ and $1 < p < 2$.
\end{proposition}

\begin{proof}
Square the left side and expand
  \begin{align*}
    \int 
    f_{z_1} g_{z_2} \overline{f_{z_3} g_{z_4}} K_N(z_1, z_2,
    z_3, z_4) \, dz_1 dz_2 dz_3 dz_4,
  \end{align*}
  where $K_N := K \chi_{|\xi_1 - \xi_2|\sim N, \ |\xi_3 - \xi_4| \sim
    N}$, and
  \begin{align*}
    K(z_1, z_2, z_3, z_4) = \langle U(t) \psi_{z_1} U(t)
    \psi_{z_2}, \ U(t) \psi_{z_3} U(t) \psi_{z_4}
    \rangle_{L^2_{t,x}([-T_0, T_0] \times \mathbf{R}^d)}.
  \end{align*}
  The estimate would follow if we could show that
  \begin{equation}
  \label{e:weighted_L2}
  N^{-\alpha}\langle z_1 - z_2 \rangle^\theta \langle z_3 - z_4
  \rangle^\theta |K_N(\vec{z})| \text{ is a bounded operator on } L^2_{z_1,
  z_2} \text{ for some } \theta > 0,
\end{equation}
as Young's inequality would then imply
\begin{align*}
\Bigl\| \int u_z \, dz \Bigr\|_{L^4}^2 &\lesssim \left( \int |f_{z_1}
  g_{z_2}|^2 \langle z_1 - z_2 \rangle^{-2\theta} \, dz_1 dz_2
  \right)^{1/2} \left( \int |f_{z_3} g_{z_4}|^2 \langle z_3 - z_4
                                         \rangle^{-2\theta} dz_3 dz_4 \right)^{1/2}\\
  &\lesssim \sup_{z} |f_z|^{2/p'} \sup_{z} |g_z|^{2/p'}
    \|f\|_{L^2}^{2/p} \|g\|_{L^2}^{2/p} \quad \text{for
    some} \quad 1 < p < 2.
\end{align*}

In view of the crude bound
$|K(\vec{z})| \lesssim \min_{j,k} \langle z_j - z_k \rangle^{-1}$,
which follows simply from the spacetime supports of the
wavepackets, ~\eqref{e:weighted_L2} would follow
from
\begin{lemma}
    \label{l:dyadic_L2}
    The localized kernel $ K_N$ satisfies
  \begin{align*}
    \| |K_N|^{1-\delta}\|_{L^2_{z_1 z_2} \to L^2_{z_3 z_4}} \lesssim N^{\alpha},
  \end{align*}
  where $\alpha$ is a constant depending only on the dimension.
\end{lemma}

\begin{proof}[Proof of Lemma~\ref{l:dyadic_L2}]
  In view of the unit scale spatial localization of the wavepackets
  and the propagation estimates~\eqref{e:integrated_bichar}, we
  may further truncate the kernel to the phase space region
  \[
   R = \{ |x_1 - x_2| \le 4|\xi_1 - \xi_2|, \quad |x_3 -
    x_4| \le 4|\xi_3 - \xi_4| \}.
  \]
  For instance, if $|x_1^s - x_2^s| \ge 4 |\xi_1^s - \xi_2^s|$ and
  $|t-s| \le T_0$ with the parameter $\eta$ in~\eqref{e:time_smallness}
  chosen sufficiently small,
  \begin{align*}
    |x_1^t - x_2^t| &\ge (1 - |t-s|^2 \| \partial^2_x V\|_{L^\infty} e^{|t-s|^2 \|
                          \partial^2_x V \|_{L^\infty}}) |x_1^s -
                          x_2^s| \\
    &- (|t-s| + |t-s|^3 \| \partial^2_x
                          V\|_{L^\infty} e^{|t-s|^2 \| \partial^2_x
      V\|_{L^\infty}}) |\xi_1^s - \xi_2^s|\\
    &\ge \frac{1}{2} |x_1^s - x_2^s| - \frac{3}{2}|t-s| |\xi_1^s -
      \xi_2^s|\\
    &\ge \frac{1}{8} |x_1^s - x_2^s|,
  \end{align*}
  therefore
  $|K_N (1-\chi_R)| \lesssim_M \langle x_1 - x_2 \rangle^{-M} \langle
  x_3 - x_4 \rangle^{-M} N^{-M}$ for any $M > 0$.  Thus it suffices to
  prove that
  \begin{align*}
    \| K_N \chi_{R}\|_{L^2 \to L^2} \lesssim N^{\alpha}.
  \end{align*}

  An estimate of this flavor was proved in the 1d
  case~\cite{mkv_inv_str_1d}. We shall argue similarly, but the proof
  is somewhat simpler since we aim for a cruder bound at this stage,
  completely ignoring temporal oscillations, and
  defer the more delicate analysis to the bilinear $L^p$ estimate.

  Partition the $4$-particle phase space $(T^* \mathbf{R}^d)^4$ according
  to the degree of physical interaction between the particles. Let
  \begin{align*}
    E_0 &= \{ \vec{z} \in (T^* \mathbf{R}^d)^4 : \min_{|t| \le T_0}
          \max_{j, k} |x_j^t - x_k^t| \le 1\},\\
    E_k &= \{ \vec{z} \in (T^* \mathbf{R}^d)^4 : 2^{k-1} < \min_{|t| \le T_0}
          \max_{j, k} |x_j^t - x_k^t| \le 2^k\},
  \end{align*}
  and decompose the kernel $K_N = \sum_{k\ge 0} K_N \chi_{E_k}$.  Then
  we have the following pointwise bound
  \begin{gather}
    \label{e:kernel_bound}
    |K(\vec{z})| \lesssim_M 2^{-kM}   \frac{ \langle \xi_1^{t(\vec{z})} +
      \xi_2^{ t(\vec{z})} - \xi_3^{ t(\vec{z})} - \xi_4^{ t(\vec{z})}
      \rangle ^{-M} }{ \langle | \xi_1^{t(\vec{z})} - \xi_2^{
        t(\vec{z} )} | + |\xi_3^{ t(\vec{z})} - \xi_4^{ t(\vec{z}
        )}| \rangle}, \quad \vec{z} \in E_k,
  \end{gather}
  where $t(\vec{z})$ is a time minimizing the ``mutual distance''
  $\max_{i,j} |x_i^t - x_j^t|$. Further, the additional localization
  to $R$ implies, by the estimates~\eqref{e:integrated_bichar}, that
  \begin{align*}
    | \xi_1^t - \xi_2^t - (\xi_1 - \xi_2)| &\lesssim \frac{1}{10} |\xi_1
                                             - \xi_2|\\
    | \xi_3^t - \xi_4^t - (\xi_3 - \xi_4)| &\lesssim \frac{1}{10} |\xi_3
                                             - \xi_4|
  \end{align*}
  for all $|t| \le T_0$. In particular
  $|\xi_1^{t(\vec{z})} - \xi_2^{t(\vec{z})}| \sim |\xi_3^{t(\vec{z})}
  - \xi_4^{ t(\vec{z})}| \sim N$; thus, while the $\xi_j^t$ may vary
  rapidly with time if $x_j^t$ are extremely far from the origin, the
  \emph{relative frequencies} retain the same order of magnitude.
  
  Assuming the bound~\eqref{e:kernel_bound} for the moment, we apply
  Schur's test to complete the proof of Lemma~\ref{l:dyadic_L2}.  Fix
  $(z_3, z_4)$ belonging to the projection
  $E_k \to T^*\mathbf{R}^d_{z_3} \times T^* \mathbf{R}^d_{z_4}$, define
\begin{align*}
E_k(z_3, z_4) = \{(z_1, z_2) : (z_1, z_2, z_3, z_4) \in E_k\},
\end{align*}
and let $t_1$ be the time minimizing
$|x_3^{t_1} - x_4^{t_1}| \le 2^k$. For any
$(z_1, z_2) \in E_k(z_3, z_4)$, the mutual distance $\max_{j, k}
|x_j^t - x_k^t|$ between
$x_1^t, x_2^t, x_3^t, x_4^t$ is minimized in the time window
\begin{align*}
	I = \{t : |t-t_1| \lesssim \min \bigl(1, \frac{2^k}{ |\xi_3 - \xi_4| }
\bigr) \},
\end{align*}
as for all other times we have $|x_3^t - x_4^t| \gg 2^k$ (Corollary~\ref{c:once_collision}).

We estimate the size of the level sets of $|K|$. For a
momentum $\xi \in \mathbf{R}^d$, let
$Q_{\xi} = (0, \xi) + [-1,1]^d \times [-1,1]^d \subset T^*\mathbf{R}^d$
denote the unit phase space box centered at $(0, \xi)$, and write
$\Phi^t = \Phi(t, 0)$ for the propagator on classical phase space
relative to time $0$ for the
Hamiltonian $h(x,\xi) = \tfrac{1}{2}|\xi|^2 + V(t, x)$. For $\mu_1, \mu_2 \in 
\mathbf{R}^d$,
define
\begin{align*}
  Z_{\mu_1, \mu_2} = \bigcup_{t \in I} (\Phi^t \otimes \Phi^t)^{-1}
  \Bigl( \frac{ z_3^t + z_4^t}{2} + 2^k Q_{\mu_1} \Bigr) \times \Bigl(
  \frac{ z_3^t + z_4^t }{2} + 2^k Q_{\mu_2} \Bigr).
\end{align*}
This set is depicted schematically in Figure~\ref{f:fig} when $k =
0$, and corresponds to the pairs of wave packets
$(z_1, z_2) \in E_m(z_3, z_4)$ with momenta $(\mu_1, \mu_2)$ relative
to the wavepackets $(z_3, z_4)$ at the ``collision time''
$t(\vec{z})$. 
		
\begin{figure}
  \includegraphics[scale=1]{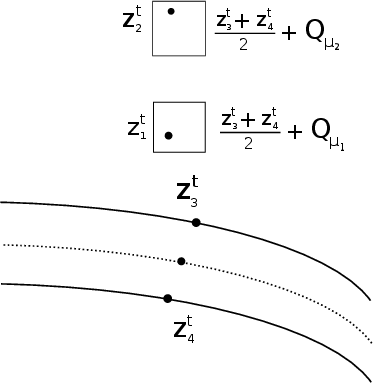}
  \caption{$Z_{\mu_1, \mu_2}$ comprises all $(z_1, z_2)$ such that 
    $z_1^t$ and $z_2^t$ belong to the depicted phase 
    space boxes for $t$ in the interval $\textrm{I}$.}
  \label{f:fig}
\end{figure}

We note that
$E_k(z_3, z_4) \subset \bigcup_{\mu_1, \mu_2 \in \mathbf{Z}^d} Z_{\mu_1,
  \mu_2}$, and recall the following estimate from the 1d paper, whose
proof we reproduce below for convenience:
\begin{lemma}
  \label{l:preimage_measure}
\begin{align}
  \label{e:preimage_measure}
  |Z_{\mu_1, \mu_2}| \lesssim 2^{4dk} \max(1, |\mu_1,|, |\mu_2|)|I|.
\end{align}
\end{lemma}

\begin{proof} Without loss assume $|\mu_1| \ge
|\mu_2|$. Partition the interval $I$ into subintervals of length
$|\mu_1|^{-1}$ if $\mu_1\neq 0$ and in subintervals of length $1$ if
$\mu_1=0$. For each $t'$ in the partition, Lemma~\ref{l:technical_lma}
implies that for some constant $C > 0$ we have
			\[
			\begin{split} \bigcup_{|t-t'| \le
\min(1,|\mu_1|^{-1})} \Phi(t)^{-1} \Bigl( \frac{z_3^t+z_4^t}{2} +
2^{k}Q_{\mu_1} \Bigr) &\subset \Phi(t')^{-1} \Bigl ( \frac{z_3^{t'} +
z_4^{t'}}{2} + C2^{k} Q_{\mu_1} \Bigr)\\ \bigcup_{|t-t'| \le
\min(1,|\mu_1|^{-1})} \Phi(t)^{-1} \Bigl( \frac{z_3^t+z_4^t}{2} +
2^{k}Q_{\mu_2} \Bigr) &\subset \Phi(t')^{-1} \Bigl ( \frac{z_3^{t'} +
z_4^{t'}}{2} + C2^{k} Q_{\mu_2} \Bigr),
			\end{split}
			\] and so
			\[
			\begin{split} &\bigcup_{|t-t'| \le
\min(1,|\mu_1|^{-1})} (\Phi(t) \otimes \Phi(t))^{-1} \Bigl(
\frac{z_3^t + z_4^t}{2} + 2^{k} Q_{\mu_1} \Bigr) \times \Bigl( \frac{z_3^t
+ z_4^t}{2} + 2^{k} Q_{\mu_2} \Bigr)\\ &\subset (\Phi(t') \otimes
\Phi(t'))^{-1} \Bigl ( \frac{z_3^{t'} + z_4^{t'}}{2} + C2^{k}
Q_{\mu_1} \Bigr) \times \Bigl ( \frac{z_3^{t'} + z_4^{t'}}{2} + C2^{k}
Q_{\mu_2} \Bigr).
			\end{split}
			\] By Liouville's theorem, the right side has
measure $O(2^{4dk})$ in $(T^* \mathbf{R}^d)^2$. The claim follows by summing
over the partition.
		\end{proof}

For each $(z_1, z_2) \in E_k(z_3, z_4) \cap Z_{\mu_1, \mu_2}$,  we
have by definition $z_j^{t(\vec{z})} \in \tfrac{z_3^{t(\vec{z})} +
  z_4^{t(\vec{z})}}{2} + 2^k Q_{\mu_j}$, thus
\begin{gather*}
  \xi_1^{t(\vec{z})} + \xi_2^{t(\vec{z})} - \xi_3^{ t(\vec{z})}  -
  \xi_4^{t(\vec{z})} = \mu_1 + \mu_2 + O(2^k)\\
  \xi_1^{t(\vec{z})} - \xi_2^{t(\vec{z})} = \mu_1 - \mu_2 + O(2^k)
\end{gather*}
Hence when $(z_1, z_2) \in Z_{\mu_1 , \mu_2}$, for any $M$ we have
\begin{align}
  \label{e:schur_kernel_bound}
  |K(\vec{z})| \lesssim_M 2^{-Mk} \frac{ \langle \mu_1 +
  \mu_2\rangle^{-M} }{\langle |\mu_1 - \mu_2| + |\xi_3^{t(\vec{z})} -
  \xi_4^{t(\vec{z})}| \rangle}.
\end{align}

To apply Schur's test, we combine the
estimates~\eqref{e:preimage_measure}, \eqref{e:schur_kernel_bound},
and evaluate
\begin{align*}
  \int |K_N (z_1, z_2, z_3, z_4)|^{1-\delta} \chi_{E_k}(\vec{z}) \,
  dz_1 dz_2 &\le \sum_{\mu_1, \mu_2 \in \mathbf{Z}^d}\int_{Z_{\mu_1,
              \mu_2}} |K_N^{1-\delta} \chi_{E_k} \, dz_2 dz_2\\
            &\lesssim_M 2^{-Mk} \sum_{|\mu_1 - \mu_2| \lesssim N +
              2^k} 2^{-Mk}  \langle \mu_1 + \mu_2 \rangle^{-M}\\
  &\lesssim N^{d} 2^{-(M-d)k}.
\end{align*}
For fixed $z_1, z_2$, the integral over $z_3$ and $z_4$ is estimated
the same way. This concludes the proof of Lemma~\ref{l:dyadic_L2},
modulo some remarks on the crucial pointwise
bound~\eqref{e:kernel_bound}.

To obtain that estimate, we use Lemma~\ref{l:galilei} to write
\begin{align*}
  K(\vec{z}) &= \int e^{i\Phi} \prod_{j=1}^4 U^{z_j} (t) \psi
               (x-x_j^t), \, dx dt,\\
  \Phi(t, x; \vec{z}) &= \sum_j \sigma_j \Bigl[\langle x-x_j^t,
                        \xi_j^t \rangle + \phi(t, x_0, \xi_0)\Bigr]
\end{align*}
where $\sigma = (+, +, -, -)$, and we denote
$\prod_j c_j := c_1 c_2 \overline{c_3} \overline{c_4}$. 

It is convenient to partition the integral further, writing
\[U^{\vec{z}_j}(t) \psi (x-x_j^t) =
\sum_{\ell_j \ge 0} U^{\vec{z}_j}(t) \psi (x-x_j^t) \theta_{\ell_j}(x-x_j^t)\bigr),
\]
where $\sum_{\ell \ge 0} \theta_{\ell}$ is a partition of unity with
$\theta_\ell$ supported on the dyadic annulus of radius
$\sim 2^{\ell}$.  For $\vec{z} \in E_k$, only the terms
\begin{align*}
  K_{\vec{\ell}}(\vec{z}) := \int e^{i\Phi} \prod_{j=1}^4 U^{z_j} (t) \psi
  (x-x_j^t) \theta_{\ell_j} (x-x_j^t) \, dx dt
\end{align*}
with $\ell^* := \max_j \ell_j \gtrsim k$ will be nonzero.

By Lemma~\ref{l:bichar}, the integral is supported on the spacetime
region
\begin{align*}
  \{(t,x) : |t - t(\vec{z})| \lesssim \min \Bigl(1, \frac{ 2^{\ell^*} }{ 
  \max_{i, j}
  |\xi_i^{t(\vec{z})} - \xi_k^{t(\vec{z})}|} \Bigr) \quad \text{and} \quad |x-x_j^t| \lesssim 2^{\ell_j}\},
\end{align*}
and for all such $t$ we have
\begin{align*}
  |x_j^t - x_k^t| \lesssim 2^{\ell^*}, \quad | \xi_j^t - \xi_k^t - (\xi_j^{t(\vec{z})}
  - \xi_k^{t(\vec{z})})| \lesssim 2^{\ell^*}.
\end{align*}
Integrating by parts in $x$, we may produce as many factors of
$|\xi_1^t + \xi_2^t - \xi_3^t - \xi_4^t |^{-1}$ as desired and freeze
$t = t(\vec{z})$ to obtain
\begin{align*}
  |K_{\vec{\ell}} (\vec{z})| \lesssim_M 2^{-\ell^* M} \frac{ \langle
  \xi_1^{t(\vec{z})} + \xi_2^{t(\vec{z})} - \xi_3^{t(\vec{z})} - \xi_4^{t(\vec{z})} \rangle^{-M}}{ \langle |\xi_1^{t(\vec{z})} - \xi_2^{t(\vec{z})}|
  + |\xi_3^{t(\vec{z})} - \xi_4^{t(\vec{z})}| \rangle} \quad \text{for any} \quad M \ge 0,
\end{align*}
and the bound~\eqref{e:kernel_bound} follows upon summing over $\vec{\ell}$.
\end{proof}
This completes the proof of Proposition~\ref{p:dy_refined_L2}.
\end{proof}

\section{Proof of Theorem~\ref{t:conditional_inv_str}}

We prove Proposition~\ref{p:reduced_thm1} and hence
Theorem~\ref{t:conditional_inv_str}.  Begin with a Whitney
decomposition of
\begin{align*}
(\mathbf{R}^d \times \mathbf{R}^d) \setminus \{ (\xi, \xi) : \xi \in 
\mathbf{R}^d\}
  = \bigcup_{N \in 2^\mathbf{Z}} \bigcup_{Q \in \mathcal{Q}_N} Q,
\end{align*}
where $\mathcal{Q}_N$ is the set of dyadic cubes in
$\mathbf{R}^d \times \mathbf{R}^d$ with diameter $\sim N$ and distance
$\sim N$ to the diagonal. 
For each $Q \in \mathcal{Q}_N$, its
characteristic function factorizes
$\chi_{N}^{Q} (\xi_1, \xi_2) = \chi_{N}^{Q,1} (\xi_2) \chi^{Q, 2}_N
(\xi_2)$, where $\chi_N^{Q,j}$ are characteristic functions of
$d$-dimensional cubes of width $N$. Then we can decompose
\begin{align*}
  1(\xi_1, \xi_2) = \chi_0(\xi_1, \xi_2) + \sum_{N \ge 1} \sum_{Q \in 
  \mathcal{Q}_N}
  \chi^{Q,1}_N(\xi_1) \chi^{Q,2}_N(\xi_2),
\end{align*}
where $\chi_0(\xi_1, \xi_2)$ is supported on the set
$|\xi_1 - \xi_2| \lesssim 1$.

Now suppose $u$ and $v$ are linear solutions with initial data
$f = \int f_z \psi_z \, dz$ and $g = \int g_z \psi_z \, dz$, respectively, where
$f_z = \langle f, \psi_z \rangle$ and $g_z = \langle g, \psi_z\rangle $. 
Writing $u_z = f_zU(t) \psi_z$, $v_z = g_z U(t) \psi_z$, we deduce as a 
consequence of Hypothesis~\ref{h:bilinear_Lp} that
\begin{align}
  \label{e:bilinear_Lp_wp}
  \Bigl\| \sum_{Q \in \mathcal{Q}_N} \int_{Q}  u_{z_1} v_{z_2} \, dz_1 dz_2
  \Bigr\|_{L^q([-T_0, T_0] \times \mathbf{R}^d)} \lesssim N^{-\delta} \| f_z 
  \|_{L^2_z} \| g_z\|_{L^2_z}.
\end{align}
for each $N \ge 1$. Indeed, for each cube $Q$ the integral has a
product structure
\begin{align*}
\int_{Q} u_{z_1} v_{z_2}
dz_1 dz_2 &= \Bigl(\int u_{z_1} \chi^{Q,1}_N(\xi_1) \, dx_1
d\xi_1\Bigr) \Bigl( \int v_{z_2} \chi^{Q, 2}_N(\xi_2) \, dx_2 d\xi_2
            \Bigr)\\
  &= U(t) \Bigl[\int f_{z_1} \chi_N^{Q,1}(\xi_1) \psi_{z_1} \, dx_1
d\xi_1\Bigr] U(t) \Bigl[\int g_{z_2} \chi_N^{Q,1}(\xi_2) \psi_{z_2} \, dx_2
d\xi_2 \Bigr].
\end{align*}
By the rapid decay of the wavepackets, we may harmlessly insert
frequency cutoffs $\tilde{\chi}_{N}^{Q, j}(D)$, where
$\tilde{\chi}_N^{Q,j}$ are slightly fattened versions of
$\chi_N^{Q, j}$ and still have supports separated by distance $\sim N$,
and apply Hypothesis~\ref{h:bilinear_Lp} to estimate
\begin{align*}
  \Bigl\| \int_Q u_{z_1} v_{z_2} \, dz_1 dz_2 \Bigr\|_{L^q} &\lesssim
  N^{-\delta} \Bigl\| \int f_{z_1} \chi^{Q,1}_N (\xi_1) \, dx_1 d\xi_1
  \Bigr\|_{L^2(\mathbf{R}^d)}  \Bigl\| \int g_{z_2} \chi^{Q,2}_N (\xi_2)
                                                              \, dx_2
                                                              d\xi_2
                                                              \Bigr\|_{L^2(\mathbf{R}^d)}\\
  &\lesssim N^{-\delta} \| f_{z} \chi_N^{Q,1}(\xi)\|_{L^2_z} \| g_z
    \chi_N^{Q,2} \chi(\xi) \|_{L^2_z}.
\end{align*}
The left side of~\eqref{e:bilinear_Lp_wp} is therefore bounded by
\begin{align*}
\sum_{Q \in \mathcal{Q}_N} N^{-\delta} \| f_{z} \chi_N^{Q,1}(\xi)\|_{L^2_z} \| 
g_z
                 \chi_N^{Q,2} \chi(\xi) \|_{L^2_z} &\le N^{-\delta}
  \Bigl( \sum_{Q \in \mathcal{Q}_N} \| f_z \chi_N^{Q, 1}(\xi)\|_{L^2_z}^2 
  \bigr)^{1/2}
  \Bigl( \sum_{Q \in \mathcal{Q}_N} \| g_z \chi_N^{Q, 2}(\xi)\|_{L^2_z}^2
  \bigr)^{1/2}\\
  &\lesssim N^{-\delta} \|f_z\|_{L^2_z} \|g_z\|_{L^2_z},
\end{align*}
as claimed.

Now decompose
the product
\begin{align*}
  uv &= \int u_{z_1} v_{z_2} \chi_0(\xi_1, \xi_2) \, dz_1 dz_2 +
       \sum_{N \ge 1} \sum_{Q \in \mathcal{Q}_N} \int_{Q}
       u_{z_1} v_{z_2} \, dz_1 dz_2,
\end{align*}
and estimate each group of terms in $L^q$ for $q$ between $p$ and
$2$. For the sum over $\mathcal{Q}_N$ we interpolate between the $L^p$
and $L^2$ bounds. Writing
$\tfrac{1}{q} = \tfrac{1-\theta}{p} + \tfrac{\theta}{2}$, we have
\begin{align*}
  \Bigl\| \sum_{Q \in \mathcal{Q}_N} \int_Q u_{z_1} v_{z_2} \, dz_1 dz_2
  \Bigr\|_{L^q} &\le \Bigl\| \sum_{Q \in \mathcal{Q}_N} \int_Q u_{z_1} v_{z_2} 
  \, dz_1 dz_2
  \Bigr\|_{L^p}^{1-\theta} \Bigl\|
                  \sum_{Q \in \mathcal{Q}_N} \int_Q u_{z_1} v_{z_2} \, dz_1 dz_2
  \Bigr\|_{L^2}^\theta\\
  &\lesssim N^{-\delta (1-\theta)+\alpha \theta} \bigl[ (\sup_z |\langle f, \psi_z|)^{1/p'}
    \sup_{z}  | \langle g, \psi_{z}|)^{1/p'} \bigr]^{\theta} (\| f\|_{L^2_x} \|
    g\|_{L^2_x})^{1-\theta + \frac{1\theta}{p}}
\end{align*}
and for $q$ sufficiently close to $p$ (hence $\theta$ sufficiently
small) the exponent of $N$ is negative.

For the ``near-diagonal'' sum, we interpolate between $L^1$
and $L^2$. For the $L^1$ bound we simply use Minkowski's inequality
and the estimate $\|U(t)\psi_{z_1} U(t) \psi_{z_2}\|_{L^1} \lesssim_N
\langle x_1 -x_2\rangle^{-N}$ when $|\xi_1-\xi_2| \le 1$ to obtain
\begin{align*}
  \Bigl \| \int u_{z_1} v_{z_2} \chi_0(\xi_1, \xi_2) dx_1 dx_2 d\xi_1
  d\xi_2 \Bigr \|_{L^1_x} &\lesssim \int |f_{z_1} g_{z_2}| \langle x_1 - x_2
                            \rangle^{-N} \chi_0(\xi_1, \xi_2) \, dz_1 dz_2\\
  &\lesssim \| f_{z} \|_{L^2_z} \| g_{z} \|_{L^2_z},
\end{align*}
which when combined with Proposition~\ref{p:dy_refined_L2} yields
\begin{align*}
  \Bigl\| \int u_{z_1} v_{z_2} \chi_0(\xi_1, \xi_2) \, dz_1
  dz_2\Bigr\|_{L^q} &\lesssim \Bigl\| \int u_{z_1} v_{z_2} \chi_0(\xi_1, \xi_2) \, dz_1
  dz_2\Bigr\|_{L^1}^{1-{\theta'}} \Bigl\| \int u_{z_1} v_{z_2} \chi_0(\xi_1, \xi_2) \, dz_1
  dz_2\Bigr\|_{L^2} ^{\theta'}\\
  &\lesssim \bigl[(\sup_{z} |f_z|)^{1/p'} (\sup_{z}
    |g_z|)^{1/p'}\bigr]^{{\theta'}} (\|f_z\|_{L^2_z} \|g_z\|_{L^2_z})^{1
    - {\theta'} + \frac{\theta'}{p}}\\
  &\lesssim (\sup_z |\langle f, \psi_z\rangle| \sup_{z} |\langle g,
    \psi_{z} \rangle|)^{\theta/p'} (\|f\|_{L^2}\|g\|_{L^2})^{1 -
    \theta + \frac{\theta}{p}}
\end{align*}
for some $1 < p < 2$, where
$\tfrac{1}{q} = 1-{\theta'} + \tfrac{{\theta'}}{2}$.

Summing in $N$, we conclude that
\begin{align*}
  \|uv\|_{L^{q}} \lesssim \bigl[(\sup_{z} |\langle f, \psi_z|)^{1/p'} (\sup_{z}
    |\langle g, \psi_z \rangle|)^{1/p'}\bigr]^{{\theta}} (\|f\|_{L^2_x} \|g\|_{L^2_x})^{1
    -\frac{\theta}{p'}}
\end{align*}
for some ${\theta} = {\theta}(p) \in (1, \tfrac{d+2}{d})$. Taking
$u = v$ we obtain Proposition~\ref{p:reduced_thm1}.

\section{The restriction-type estimate}
\label{s:restriction}

This purpose of this section is to prove Theorem~\ref{t:bilinear_Lp}.

We shall systematically use the following notation.  For $N \ge 1$
and a potential $V$, we consider the rescaled potentials
\[V_N(t, x) : = N^{-2} V(N^{-2}t, N^{-1}x).\] Let $U(t, s)$ and
$U_N(t, s)$ denote the propagators for the corresponding
Schr\"{o}dinger operators $H(t) : = -\tfrac{1}{2}\Delta + V$ and
$H_N(t) := -\frac{1}{2} \Delta + V_N$. We will often use the letter $U$
to write the propagators for different potentials $V \in \mathcal{V}$;
this ambiguity will not cause any serious issue, however, since all
the estimates we shall need are valid uniformly over
$\mathcal{V}$. Further, due to the time translation invariance of our
assumptions we shall usually just consider the propagator from time
$0$ and write $U(t) : = U(t, 0)$, $U_N(t) := U_N(t, 0)$.

In the sequel, the letter $C$ will denote a constant, depending only
on the dimension $d$, which may change from line to line.

\subsection{Preliminary reductions}

The hypotheses of Theorem~\ref{t:bilinear_Lp} are invariant under
various transformations of $u$ and $v$.
\begin{itemize}
\item Galilei boosts $u(0) \mapsto \pi(z_0) u(0), \ u \mapsto
  \pi(z_0^t)u^{z_0}$, where $u^{z_0}$ satisfies $(D_t - \Delta +
  V^{z_0}) u^{z_0} = 0, \ u^{z_0}(0) = u(0)$.
\item Spatial rotations: for an orthogonal matrix
  $g$, $(g \cdot u)(t, x) := u(t, g^{-1}\cdot x)$ satisfies
  \begin{align*}
    [D_t (g \cdot u) - \Delta+ (g \cdot V)] (g \cdot u) = 0.
  \end{align*}
\item Rescaling
  $u \mapsto u_\lambda = \lambda^{-\frac{d}{2}} u(\lambda^{-2} t, \lambda^{-1} x)$ for
  $\lambda > 1$. Then $u_\lambda$ satisfies $(D_t -\Delta + V_\lambda)u_\lambda = 0$ with a
  smoother potential $V_\lambda(t,x) = \lambda^{-2} V(\lambda^{-2} t, \lambda^{-1} x)$.
\end{itemize}

We may and shall assume hereafter that $V$ vanishes to second order at
  $x=0$, that is, $V(t, 0) = 0$ and $\partial_x V (t, 0) = 0$ for all
  $t$. Indeed let $z_0^t = (x_0^t, \xi_0^t)$ be the bicharacteristic with $(x_0, \xi_0) =
  (0, 0)$. Then by Lemma~\ref{l:galilei},
  \begin{align*}
   \| U(t) f U(t) g \|_{L^{\frac{d+3}{d+1}} } = \| (\pi(z_0^t)
    U^{z_0}(t) f) (\pi(z_0^t) U^{z_0}(t) g)\|_{L^{\frac{d+3}{d+1}}} =
    \| U^{z_0}(t) f U^{z_0} (t) g\|_{L^{\frac{d+3}{d+1}}},
  \end{align*}
and the potential $V^{z_0}(t,x) = V(t, x_0^t + x) - V(t, x_0^t) - x
\partial_x V (t, x_0^t)$ vanishes to second order at $x = 0$.

Theorem~\ref{t:bilinear_Lp} is equivalent by rescaling to
\begin{theorem}
  \label{t:bilinear_Lp_rescaled}
  Given $S_1, S_2 \subset \mathbf{R}^d_{\xi}$ with $\operatorname{diam}(S_j) 
  \le 1$
  and $c^{-1} \ge \operatorname{dist}(S_1, S_2) \ge c$ for some $0 < c < 1$,
  there exists a constant $\eta = \eta(c) > 0$ such that if $V \in
  \mathcal{V}$ and $\tau_0>0$ satisfies
  \begin{align}
    \label{e:Vsmallness}
    (\tau_0 + \tau_0^2) \| \partial^2_x V\|_{L^\infty_{t,x}} < \eta,
  \end{align}
  then for any $f, g \in L^2(\mathbf{R}^d)$ with
  $\operatorname{supp}(\hat{f}) \subset S_1$ and
  $\operatorname{supp}(\hat{g}) \subset S_2$, the corresponding
  Schr\"{o}dinger solutions $u_N = U_N(t) f$ and $v_N = U_N(t) g$
  satisfy the estimate
  \begin{align}
    \label{e:bilinear_Lp_rescaled}
    \| u_N v_N\|_{L^q ( [-\tau_0 N^2, \tau_0 N^2] \times \mathbf{R}^d)}
    \lesssim_{\varepsilon} N^\varepsilon \|f\|_{L^2} \|g\|_{L^2} \quad
    \text{for all} \quad \frac{d+3}{d+1} \le q < \frac{d+2}{d}
  \end{align}
  for any $\varepsilon > 0$ and $N \ge 1$.
\end{theorem}

\ifdraft
[DRAFT]
\begin{align}
  \label{e:integrated_bichar_V}
  x_1^t - x_2^t = (1 + O( |t-s|^2\|\partial^2_xV\|) )(x_1^s - x_2^s  ) +
  (I + O(|t-s|^2 \| \partial^2_xV\|) )(t-s) (\xi_1^s - \xi_2^s).
\end{align}

Observe, the definition of the
rescaled potential $V_N$ implies
\[
  |t| \le \tau_0 N^2 \Rightarrow t^2 \| \partial^2_x V_N\|_{L^\infty} < \eta.
\]
\fi

In fact it suffices to take $S_1$ and $S_2$ of the form
\begin{align}
  \label{e:small_fourier_supports}
  S_1 = \{\xi: |\xi - \frac{c}{2}e_1| \le \frac{c}{100}\}, \  S_2 =
  \{\xi: |\xi + \frac{c}{2}e_1| \le \frac{c}{100}\}.
\end{align}
General $S_j$ can be reduced to this case by decomposing
$\hat{f} = \sum_j \widehat{f_j}$ and $\hat{g} = \sum_k \widehat{g_k}$
into pieces supported in small balls and applying an appropriate
Galilei boost and rotation for each pair $(f_j, g_k)$ and possibly
also a rescaling to bring the Fourier supports closer, which only
reduces $\|\partial^2_xV\|_{L^\infty}$. Henceforth we shall
assume~\eqref{e:small_fourier_supports}.

\subsection{General remarks}
We use the induction on scales method pioneered by Wolff for the
cone~\cite{Wolff2001} and adapted by Tao to the
paraboloid~\cite{Tao2003}. Our proof is modeled closely on Tao's
treatment of the $V=0$ case, and the reader may find it helpful to
read the following exposition in parallel with~\cite{Tao2003}. The
main differences are as follows:

\begin{itemize}
\item The induction scheme (section~\ref{ss:induction}) is complicated
  by the fact that frequency is not conserved, so one cannot directly
  apply an induction hypothesis which involves assumptions on the
  frequency supports at time $0$ to a spacetime ball at a later time.

\item The low regularity of $V$ in time makes the bilinear $L^2$
  estimate (section~\ref{ss:restriction_L2_bound}) more delicate and
  we obtain weaker decay from temporal oscillations.

\item In the final Kakeya-type estimate, the tubes in the key combinatorial
  lemma (Lemma~\ref{l:tube_comb}, the analogue of Lemma 8.1 in Tao)
  are curved. Also, we need to be slightly more precise to compensate
  for the weaker decay in the $L^2$ bound.
    
\end{itemize}

\subsection{Discrete wavepacket decomposition}
\label{ss:discrete_wp_decomp}

While the first part of this paper employed continuous wavepacket
transforms, the following discrete decomposition, taken essentially
from Tao~\cite{Tao2003}, is more conventional in restriction theory
and convenient for the combinatorial arguments involved. To each
$z_0 = (x_0, \xi_0)$ in classical phase space with bicharacteristic
$\gamma_{z_0}(t)= (x_0^t, \xi_0^t)$, we associate a spacetime ``tube''
\begin{align*}
  T_{z_0} := \{(t, x) : |x-x_0^t| \le R^{1/2}, \quad |t| \le R\}.
\end{align*}
For such a tube $T$, let $z(T) = (x(T), \xi(T))$ denote the
corresponding initial point in phase space. A wavepacket $\phi$
associated to the bicharacteristic $z_0 \mapsto z_0^t$ is essentially
supported in spacetime on the tube $T_{z_0}$, and we shall often
emphasize this fact by writing $\phi_T$.

\begin{lemma}
  \label{l:wp_decomp}
  Let $u = U_N(t) f$ be a linear Schr\"{o}dinger solution with
  $\operatorname{supp}(\hat{f}) \subset S_1$. For each $1\le R\le N^2$, there
  exists a collection of tubes $\mathbf{T}$ and a decomposition
  \begin{align*}
    u = \sum_{T \in \mathbf{T}} a_T \phi_T,
  \end{align*}
  into $R \times (R^{1/2})^d$ wave packets with the following properties:
  \begin{itemize}
  \item Each $T \in \mathbf{T}$ satisfies $( x(T), \xi(T) )
    \in R^{1/2} \mathbf{Z}^d \times R^{-1/2} \mathbf{Z}^d$. 
  \item Each wavepacket $\phi_T$ is a Schr\"{o}dinger solution
    localized near the bicharacteristic $(x(T)^t, \xi(T)^t)$,
    i.e. which satisfies the pointwise bounds
    \begin{equation}
      \label{e:wp_bounds}
      \begin{split}
  |(R^{1/2} \partial_x)^k \phi_{T}(t)| &\lesssim_{k,M} \Bigl \langle
      \frac{x-x(T)^t}{R^{1/2}} \Bigr \rangle^{-M}, \\
      |(R^{-1/2}
  \partial_\xi)^k \widehat{\phi_{T}}(t)| &\lesssim_{k,M} \Bigl \langle
  \frac{\xi-\xi(T)^t}{R^{-1/2}} \Bigr\rangle^{-M} 
  \end{split} \quad \text{ for all } k, \ M\ge
  0.
\end{equation}
Moreover, $\widehat{\phi_T}[0]$ is supported in a $R^{-1/2}$
neighborhood of $\xi(T) \in S_1$.
\item The complex coefficients $a_T$ are square-summable:
  \begin{align*}
    \sum_{T} |a_T|^2 \lesssim \|f\|_{L^2}^2.
  \end{align*}
Moreover, for any subcollection of tubes $\mathbf{T}' \subset \mathbf{T}$ and
  complex numbers $a_T$, one has
  \begin{align*}
    \|\sum_{T \in \mathbf{T}'} a_T \phi_T \|_{L^2}^2 \lesssim \sum_{T \in
    \mathbf{T}'} |a_T|^2.
  \end{align*}
\end{itemize}
A similar decomposition also holds for $v = U_N(t)g$.
\end{lemma}

\begin{proof}[Proof sketch]
  We outline the main steps as this construction is fairly standard;
  consult for instance Lemma~4.1 in~\cite{Tao2003}.  Begin with
  partitions of unity $1 = \sum_{x_0 \in \mathbf{Z}^d} \eta(x - x_0)$ and
  $1= \sum_{\xi_0 \in \mathbf{Z}^d} \chi(\xi - \xi_0)$ such that $\chi$ and
  $\hat{\eta}$ are compactly supported. By rescaling and quantizing,
  we obtain a pseudo-differential partition of unity used to decompose
  the initial data
  \begin{align*}
    f = \sum_{(x_0, \xi_0) } \eta\Bigl( \frac{x-x_0}{R^{1/2}} \Bigr) \chi(
    R^{1/2} (D - \xi_0)) f.
  \end{align*}
  The propagation estimates then follow from the next lemma.
\end{proof}

\begin{lemma} If $\phi_{z_0}$ is a scale-$R$ wavepacket concentrated
  at $z_0$, and $U_N(t)$ is the propagator for $H(t) = -\tfrac{1}{2}\Delta
  + V_N$, then $U_N(t)$ is a scale-$R$ wavepacket concentrated at
  $z_0^t$ for all $|t| \le R$.
\end{lemma}

\begin{proof} By rescaling we reduce to $R = 1$ and replace $V$ by
  $V_{N/R^{1/2}}$ which also belongs to $\mathcal{V}$ since $N/R^{1/2}
  \ge 1$. Then the symbol $a = \tfrac{1}{2}|\xi|^2 + V_{N/R^{1/2}}(t,x)$
  satisfies the estimates~\eqref{e:intro_symbol_hyp}, and we can appeal
  to Lemma~\ref{l:wp_propagation}.
\end{proof}

\subsection{Localization}
\label{s:localization}

The proof of Theorem~\ref{t:bilinear_Lp_rescaled} begins with the
observation that it suffices to establish the same estimate with the
spacetime norm restricted to a box of the form
\[\Omega_N = [-N^2, N^2] \times [-AN^2, AN^2]^d.\]

\begin{theorem}
\label{t:loc_bilinear_Lp}
Assume the hypotheses and notation of
Theorem~\ref{t:bilinear_Lp_rescaled} and replace $c$ by $c/2$ and take
$\operatorname{diam}(S_j) \le 11/10$. Then there exists $A  = A(c)> 0$ such that
  \begin{align}
    \label{e:loc_bilinear_Lp}
    \| u_N v_N \|_{L^{\frac{d+3}{d+1}}( \Omega_N)}
    \lesssim_\varepsilon N^\varepsilon \| f\|_{L^2} \| g\|_{L^2} 
  \end{align}
  for any $\varepsilon > 0$.
\end{theorem}
\begin{remark}
  In the wavepacket decomposition of $u_N$ and $v_N$, the Fourier
  supports of the wavepackets are contained in a slight dilate
  $S_j + B(0, CN^{-1})$ of $S_j$. Hence at various junctures we need
  to adjust various constants to accommodate this minor enlargement of
  Fourier supports.
\end{remark}

The full theorem then follows from an approximate finite speed of
propagation argument:
\begin{lemma}
  Theorem~\ref{t:loc_bilinear_Lp} implies Theorem~\ref{t:bilinear_Lp_rescaled}.
\end{lemma}

\begin{proof}[Proof of Lemma]
  Partition physical space
  $\mathbf{R}^d = \bigcup_{j \in \mathbf{Z}^d} Q_j$ into cubes of width
  $\sim N^2$, where $Q_j$ denotes the cube with center
  $N^2j \in N^2 \mathbf{Z}^d$. Decompose $u := u_N$ and
  $v := v_N$ into $N^2 \times (N)^d$ wavepackets, and group the terms
  in the product according to their relative initial positions. Write
\begin{align*}
  u &= \sum_T a_T \phi_T = \sum_{j \in \mathbf{Z}^d} \sum_{T \in \mathbf{T}_j} 
  u_T,\\
  v &= \sum_{T'} b_{T'} \phi_{T'} = \sum_{j' \in \mathbf{Z}^d} \sum_{T' \in 
  \mathbf{T}'_{j'}} v_{T'},
\end{align*}
where $\mathbf{T}_j = \{ T \in \mathbf{T} : x(T) \in Q_j\}$ and similarly for
$\mathbf{T}_{j'}$.
Using the triangle inequality we estimate
\begin{align}
  \label{e:localization-sum}
\|  uv \|_{L^{\frac{d+3}{d+1}}}&\le   \sum_{k \ge 0} \Bigl\| \sum_{|j - j'| \sim 2^k} \sum_{T \in
       \mathbf{T}_j, \ T' \in \mathbf{T}'_{j'}} u_T v_{T'} 
       \Bigr\|_{L^{\frac{d+3}{d+1}}}.
\end{align}
For the $k$th sum, note from~\eqref{e:integrated_bichar} that if
$(x_1, \xi_1) := (x(T), \xi(T))$ and
$(x_2, \xi_2) := (x(T'), \xi(T'))$, we have
\begin{align*}
  |x_1^t - x_2^t| &\ge (1 - Ct^2 \| \partial^2_x V_N\|_{L^\infty})  |x_1 - x_2| - (|t| + C|t|^3 \|
  \partial^2_x V_N\|_{L^\infty}) |\xi_1 - \xi_2|\\
  &\ge (1 -  C\tau_0^2\| \partial^2_x V\|_{L^\infty})|x_1 - x_2| - N^2 (1 +
     C \tau_0^2 \|\partial^2_x V\|_{L^\infty}|\xi_1 - \xi_2|\\
  &\ge (1 - C\eta)|x_1 - x_2| - N^2(1 + C\eta)|\xi_1 - \xi_2|,
\end{align*}
where $C$ hides the harmless Gronwall factor.  As
$|\xi_1-\xi_2| \le c^{-1}$, there exists $k(c)$ such that if
$|x_1 - x_2| \ge 2^kN^2$ and $\eta$ is chosen small enough we obtain
$ |x_1^t - x_2^t| \gtrsim 2^k N^2$ for $k \ge k(c)$. Thus the tubes in
$\mathbf{T}_j$ and $\mathbf{T}_j'$ are separated in space by distance
$\gtrsim 2^kN^2$, and since each wavepacket $\phi_T$ decays rapidly
away from its tube $T$ in units of $N$, we have
\begin{align*}
  \|\phi_T \phi_{T'}\|_{L^{\frac{d+3}{d+1}}} &\lesssim 2^{-101dk}
                                               N^{-101d},
\end{align*}
and estimate crudely as follows:
\begin{align*}
  \Bigl\| \sum_{|j-j'| \sim 2^k}  \sum_{T \in \mathbf{T}_j,  \, T' \in \mathbf{T}'_j} u_T v_{T'}\Bigr\|_{L^{\frac{d+3}{d+1}}}
  &\lesssim   2^{-101dk} N^{-101d}  \sum_{|j-j'| \sim 2^k} \sum_{T \in \mathbf{T}_j,  \, T' \in
    \mathbf{T}'_{j'}} |a_T b_{T'}| \\
  &\lesssim 2^{-101dk} N^{-100d}  \sum_{|j-j'| \sim 2^k} \bigl (\sum_{T
    \in \mathbf{T}_j} |a_T|^2\bigr)^{\frac{1}{2}} \bigl (\sum_{T'
    \in \mathbf{T}'_{j'}} |b_{T'}|^2\bigr)^{\frac{1}{2}}\\
  &\lesssim 2^{-100dk} N^{-100d} \Bigl( \sum_{j} \sum_{T \in
    \mathbf{T}_j} |a_T|^2 \Bigr)^{\frac{1}{2}}\Bigl( \sum_{j} \sum_{T \in
    \mathbf{T}'_{j'}} |b_{T'}|^2 \Bigr)^{\frac{1}{2}}\\
  &\lesssim 2^{-{100dk}} N^{-100d}
  \|f\|_{L^2} \|g\|_{L^2}.
\end{align*}

For the ``near diagonal'' part of the sum~\eqref{e:localization-sum},
where $|j-j'| \le 2^{k(c)}$, we group the terms by their average initial
positions:
\begin{align}
  \label{e:localization}
  \Bigl \| \sum_{|j - j'| \lesssim 1} \sum_{T \in \mathbf{T}_j, T' \in
  \mathbf{T}'_{j'}} u_T v_{T'} \Bigr\|_{L^{\frac{d+3}{d+1}}} &\le \sum_{m \in 
  \mathbf{Z}^d +
                                    \mathbf{Z}^d} \sum_{|j-j'|
                                    \lesssim 1, j+j' = m} \Bigl \| \sum_{T \in
                                    \mathbf{T}_j, \, T' \in \mathbf{T}'_{j'}}
                                    u_T v_{T'}
                                                          \Bigr\|_{L^{\frac{d+3}{d+1}}}.
\end{align}
For each pair $(j, j')$, we translate the initial data by the midpoint
$x_{jj'} := \tfrac{j+j'}{2}N^2$ of $Q_j$ and $Q_{j'}$,  using
Lemma~\ref{l:galilei} to write
\[u_{T} = \pi(z_{jj'}^t) a_T \tilde{\phi}_T =: \tilde{u}_T, \
v_T = b_{T'} \pi(z_{jj'}^t) \tilde{\phi}_{T'} =: \tilde{v}_T,\] where
$z_{jj'} = (x_{jj'}, 0)$ and
\begin{align*}
  \tilde{\phi}_T(t) =  U^{(x_{jj'}, 0)}(t) \pi(-x_{jj'}, 0)
  \phi_T[0]
\end{align*}
is a wavepacket solution for the modified potential
$V^{(x_{jj'}, 0)}$. 
The norm on the right side above therefore can be written as
\begin{align*}
  \Bigl\| \sum_{T \in \widetilde{\mathbf{T}_j}, T' \in 
  \widetilde{\mathbf{T}'_j}} \tilde{u}_T
  \tilde{v}_{T'} \Bigr\|_{L^{\frac{d+3}{d+1}}},
\end{align*}
where the initial positions $x(T)$ and $x(T')$ of the tubes now belong
to the translated cubes $\tilde{Q}_j := Q_j - x_{jj'}$,
$\tilde{Q}_{j'} - x_{jj'}$, which are now distance $\lesssim
N^2$ from the origin (note however that the tubes in $\widetilde{\mathbf{T}}_j$ 
are not
simply translates of those in $\mathbf{T}_j$).

By simple bicharacteristic estimates and the wavepacket
bounds~\eqref{e:wp_bounds}, for large $A$ the norm outside
$\Omega_N := [-N^2, N^2] \times [-AN^2, AN^2]^{d}$ is negligible:
\begin{align*}
  \Bigl\| \sum_{T \in \tilde{\mathbf{T}}_j, T' \in \tilde{\mathbf{T}}'_j} 
  \tilde{u}_T
  \tilde{v}_{T'} \Bigr\|_{L^{\frac{d+3}{d+1}}( [-N^2, N^2] \times ( [-AN^2,
  AN^2]^c))} &\lesssim N^{-100d} \Bigl( \sum_{T \in \tilde{\mathbf{T}}_j}
  |a_T|^2 \Bigr)^{1/2} \Bigl( \sum_{T' \in \tilde{\mathbf{T}}'_j} |b_T|^2
  \Bigr)^{1/2}\\
  &\lesssim N^{-100d} \Bigl( \sum_{T \in \mathbf{T}_j}
  |a_T|^2 \Bigr)^{1/2} \Bigl( \sum_{T' \in \mathbf{T}'_j} |b_T|^2
  \Bigr)^{1/2}
\end{align*}
Inside $\Omega_N$ we invoke Proposition~\eqref{t:loc_bilinear_Lp}
using the fact that the $V^{(x_{jj'}, 0)}$ also satisfies the
hypothesis~\eqref{e:Vsmallness}, and that the wavepacket
decompositions of $u_N$ and $v_N$ satisfy the relaxed Fourier support
conditions in that proposition. Altogether, the right side
of~\eqref{e:localization} is bounded by
\begin{align*}
  &\sum_{m \in \mathbf{Z}^d + \mathbf{Z}^d} \sum_{|j - j'| \lesssim 1, \, j +
  j' = m} N^{\varepsilon}  \Bigl( \sum_{T \in \mathbf{T}_j}
  |a_T|^2 \Bigr )^{1/2} \Bigl( \sum_{T' \in \mathbf{T}'_j} |b_T|^2
  \Bigr)^{1/2}\\
  &\lesssim N^{\varepsilon} \sum_{m} \Bigl( \sum_{|j -
    \frac{m}{2}| \lesssim 1} \sum_{T \in \mathbf{T}_j} |a_T|^2 \Bigr)^{1/2} 
    \Bigl( \sum_{|j' -
    \frac{m}{2}| \lesssim 1} \sum_{T \in \mathbf{T}'_{j'}} |b_{T'}|^2
    \Bigr)^{1/2}\\
  &\lesssim N^\varepsilon \Bigl( \sum_{T} |a_T|^2 \Bigr)^{1/2} \Bigl(
    \sum_{T'} |b_{T'}|^2 \Bigr)^{1/2}\\
  &\lesssim N^\varepsilon \|f\|_{L^2} \|g\|_{L^2},
\end{align*}
thus recovering Theorem~\ref{t:bilinear_Lp_rescaled}.
\end{proof}

\subsection{Induction on scales}
\label{ss:induction}

Our induction scheme is set up slightly differently from Tao's to
accommodate the non-conservation of frequency support of solutions.


In this section, we explicitly display the dependence of the
propagator on the potential, and write $U^V_N(t) = U^V_N(t, 0)$ for
the propagator with potential $V_N$.

Let $\text{IH}(\alpha)$ denote the following statement:
\begin{quote}
 There exists $C_\alpha > 0$ such that for each $N \ge 1$ and for all potentials $V \in
\mathcal{V}_\eta$,
the estimate
\begin{align}
  \label{e:IH}
  \|U^V_N(t) f U^V_N(t) g\|_{L^{\frac{d+3}{d+1}}(\Omega_N)}
\le C_\alpha N^{2\alpha} \| f\|_{L^2} \|g\|_{L^2}
  \end{align} 
  holds for all $f, g \in L^2(\mathbf{R}^d)$ with $\hat{f}, \hat{g}$ supported
  in $S_1$ and $S_2$, respectively.
\end{quote}

We prove:

\textbf{Inductive Step}: If $\text{IH}(\alpha)$ holds, then
  $\text{IH}(\max\bigl( (1-\delta)\alpha, C\delta) + \varepsilon)$ holds
  for all $0 < \delta, \varepsilon \ll 1$.

By choosing $\delta$ and $\varepsilon$ sufficiently small depending
on $\alpha$, we can always arrange that $\max \bigl(
(1-\delta)\alpha,C\delta\bigr) + C\varepsilon < \alpha - c\alpha^2$ for
some absolute constant $c$, and
Theorem~\ref{t:loc_bilinear_Lp} follows.

The inductive hypothesis $\text{IH}(\alpha)$ shall be used to improve
the estimate~\eqref{e:IH} over subregions $Q_R \subset \Omega_N$ at
smaller scales $\operatorname{diam}(Q_R) \sim N^{2(1-\delta)} \ll N^2$.
\begin{proposition}
  \label{p:inductive_hyp_cor}
  Suppose $\text{IH}(\alpha)$ holds. Then for all $1 \le R \le
  N^2/16$ and all spacetime balls $Q_R \subset 2\Omega_N$ of diameter
  $R$, the estimate
  \begin{align*}
  \|U^V_N(t) f U^V_N(t) g\|_{L^{\frac{d+3}{d+1}}(Q_R)}
\le C_\alpha R^{\alpha} \| f\|_{L^2} \|g\|_{L^2}
  \end{align*} 
  holds for all $f, g \in L^2(\mathbf{R}^d)$ with $\hat{f}, \hat{g}$ supported
  in $\tilde{S}_1 := S_1+B(0, \tfrac{c}{100})$ and $\tilde{S}_2 := S_2
  + B(0, \tfrac{c}{100})$, respectively.
\end{proposition}
\begin{proof}

We begin by estimating how much the Fourier supports can shift.

\begin{lemma}
  \label{l:fourier_truncation}
  For $1 \le R \le N^2$, let $Q_R \subset 2\Omega_N$ be a spacetime
  ball with center $(t_Q, x_Q)$ and diameter $R$.  Suppose the initial
  data $f, g$ satisfy
  $\operatorname{supp}(\hat{f}) \subset \tilde{S}_1$ and
  $\operatorname{supp}(\hat{g}) \subset \tilde{S}_2$. There exist
  decompositions $u(t_Q) = f_1 + f_2$ and $v(t_Q) = g_1 + g_2$, with
  the following properties:
  \begin{itemize}
  \item $\hat{f}_1$ and $\hat{g}_1$ are supported in sets
    $S'_1, S'_2$ with
    $\operatorname{diam}(S'_j) \le \tfrac{c}{10}$ and
    $\operatorname{dist}(S'_1, S'_2) \in [\tfrac{4c}{5},
    \tfrac{5c}{4}]$.
  \item $\|f_2\|_{L^2} \lesssim N^{-100d}\|f\|_{L^2}$ and
    $\|g_2\|_{L^2} \lesssim N^{-100d} \|g\|_{L^2}$.
  \end{itemize}
\end{lemma}

\begin{proof}
  Begin  by decomposing $u = U^V_N f$ and $v =
U^V_N g$ into $N^2 \times (N)^d$ wavepackets:
\begin{align}
  \label{e:inductive_wp}
  u = \sum_{T \in \mathbf{T}_1} a_T \phi_T, \quad v = \sum_{T \in
  \mathbf{T}_2} b_T \phi_T.
\end{align}
By the spatial localization~\eqref{e:wp_bounds}, we may ignore in $u$
and $v$ the
packets whose tubes $T \in \mathbf{T}_j$ do not intersect
$2Q_N :=[-N^2, N^2] \times [-2AN^2, 2AN^2]$, as the portion of the
sum involving those terms contributes at most $O(N^{-100d}) \|f\|_{L^2}
\|g\|_{L^2}$. Thus there are $O(N^{2d})$ remaining terms.

Suppose $\phi_{T_1}$ and $\phi_{T_2}$ are wavepackets in the
decomposition for $u$.

Let $(x_1^t, \xi_1^t)$ and $(x_2^t, \xi_2^t)$ be bicharacteristics
with $|x_1|, |x_2| \le 2 AN^2$.
By ~\eqref{e:integrated_bichar}, for $|t| \le \tau_0N^2$ we have
\begin{align*}
  |\xi_1^t - \xi_2^t - (\xi_1 - \xi_2)| &\le C\tau_0 N^2 N^{-4} \|
                                          \partial^2_xV\|_{L^\infty} ( 2AN^2 + \tau_0N^2 |\xi_1 - \xi_2| ) \\
                                        &\le C (\tau_0A + \tau_0^2) \| \partial^2_x V\|_{L^\infty} \le C\eta.
\end{align*}
Therefore, recalling the definitions of $\tilde{S}_j$, we see that we have
$|\xi_1^{t_Q} - \xi_2^{t_Q}| \le \frac{c}{20} + C\eta $ if
$\xi_1, \xi_2$ both belong to $\tilde{S}_1$ or $\tilde{S}_2$, while
$|\xi_1^{t_Q} - \xi_2^{t_Q}| \in [\tfrac{9c}{10}, \tfrac{10c}{9}]$ if
$\xi_1 \in \tilde{S}_1$ and $\xi_2 \in
\tilde{S}_2$. Choosing $\eta = \eta(c)$ sufficiently small, 

Consequently, if
\begin{align}
  \label{e:freq_set_image}
  \tilde{S}_j^t := \{ \xi^t :  \xi \in \tilde{S}_j, \quad |x|
  \le AN^2\}
  \end{align}
denotes the set of frequencies of the wavepackets at time $t$, then
$  \operatorname{diam}(\tilde{S}_j^t) \le \operatorname{diam}(S_j) + C\eta$ and 
$\operatorname{dist}(\tilde{S}_1^t,
\tilde{S}_2^t) \ge \tfrac{9}{10}\operatorname{dist}(S_1, S_2)$. Now let $S'_j$
denote  $O(N^{-9/10})$ neighborhoods of $\tilde{S}_j^t$, and decompose
\begin{align*}
  u(t_Q) = f_1 + f_2, \quad v(t_Q) = g_1 + g_2,
\end{align*}
where $\widehat{f_1}$ is supported on $\tilde{S}_1$ and
$\widehat{f_2}$ on the complement, and similarly for $g_1, g_2$. For
$N$ large enough we have $\operatorname{dist}(S_1'
S_2') \in [\tfrac{4c}{5}, \tfrac{5c}{4}]$. The
estimates in the second bullet point now follow from the
rapid decay of each wavepacket from its central
frequency on the $N^{-1}$ scale (the estimates~\eqref{e:wp_bounds}
with $R = N^{2}$).
\end{proof}

The proof of the proposition concludes with several applications of
Lemma~\ref{l:galilei}.  Write
  \begin{align*}
    U(t, t_Q)f_1 = U(t,
    t_Q)\pi(x_Q,0) \pi(-x_Q, 0) f_1 = \pi(z_Q^t) U^{z_Q}(t,
    t_Q) \tilde{f}_1 = \pi(z_Q^t) \tilde{u}(t+t_q)
  \end{align*}
  where $z_Q^{t_Q} = (x_Q, 0)$.  For $|t-t_Q| \le R$ and
  $|x_Q| \le AN^2$ we have $|x_Q^t-x_Q^{t_Q}| \le 2|t-t_Q| \le 2R$
  provided that $\eta$ is sufficiently small. Therefore, denoting
  $\tilde{Q}_R = 2(Q_R - (t_Q, x_Q))$,
  \begin{align*}
    \|uv\|_{L^{\frac{d+3}{d+1}} (Q_R)} \lesssim \| \tilde{u}
    \tilde{v}\|_{L^{\frac{d+3}{d+1}} (\tilde{Q}_R)} +
    N^{-100d} \| f\|_{L^2} \|g\|_{L^2}
  \end{align*}
  It remains to consider the first term on the right side. The
  initial data $\tilde{f}_1$, $\tilde{g}_1$ for $\tilde{u}$ and
  $\tilde{v}$ have Fourier transforms supported in $S_1', S_2'$. We
  abuse notation and redenote
  \begin{align*}
    f := \tilde{f}_1, \ g := \tilde{g}_1.
  \end{align*}

 Cover $S_j' = \bigcup_{k} B_{j, k}$ by finitely overlapping
  balls of radius $\tfrac{c}{200}$. Using a subordinate partition of
  unity, we reduce to the case where
  $\supp \hat{f} \subset B_{1, k_1}$ and
  $\supp \hat{g} \subset B_{2, k_2}$. Again using
  Lemma~\ref{l:galilei}, we may assume $B_{1, k_1} = -B_{2, k_2}$ and
  that their centers lie on the $e_1$ axis.
  
  Since $2c \ge \operatorname{dist}(B_{1, k_1}, B_{2, k_2}) \ge \tfrac{c}{2}$,
  there exists some scaling factor $\lambda \in [\tfrac{1}{2}, 2]$ such that
  $\lambda^{-1} B_{j, k_j} \subset S_j$. Consider the rescalings
  \begin{align*}
    u_\lambda = U^V_{\frac{N}{\lambda}} (t)f_\lambda =
    U^{\tilde{V}}_{(2R)^{\frac{1}{2}}}(t)f_\lambda, \quad v_\lambda =
    U^V_{\frac{N}{\lambda}} (t)g_\lambda = U^{\tilde{V}}_{(2R)^{\frac{1}{2}}}(t)g_\lambda,
  \end{align*}
  where 
  \begin{align*}
    \tilde{V}(t, x) = 2R \lambda^2 N^{-2} V(2R \lambda^2 N^{-2} t,
    (2R)^{\frac{1}{2}} \lambda N^{-1} x).
  \end{align*}
  The potential $\tilde{V}$ satisfies $\| \partial^2_x \tilde{V}\|_{L^\infty}
  \le \| \partial^2_xV\|_{L^\infty}$  since
  $2R \lambda^2 N^{-2} \le 8 RN^{-2} \le \tfrac{1}{2}$, and
  $\widehat{u}_\lambda(0)$ and $\widehat{v}_\lambda(0)$ are supported
  in $S_1$ and $S_2$. Hence we can apply $\text{IH}(\alpha)$ to
  conclude that
  \begin{align*}
  \|\tilde{u} \tilde{v}\|_{L^{\frac{d+3}{d+1}} (\tilde{Q}_R)} \lesssim
    \|u_\lambda v_\lambda\|_{L^{\frac{d+3}{d+1}} (\tilde{Q}_{2R})} \le C_\alpha
    R^{\alpha} \|f_\lambda\|_{L^2} \| g_\lambda\|_{L^2}.
  \end{align*}
\end{proof}


From here on the argument hews closely to Tao's. We recall the
following notation: write
\[A \lessapprox B\] if $A \lesssim_\varepsilon N^\varepsilon B$ for all
$N \gg 1 $ and for all $\varepsilon > 0$.

To reiterate, we want to prove
\begin{align}
  \label{e:reduced_Lp_bilinear}
  \| U^V_N f U^V_N g\|_{L^{\frac{d+3}{d+1}} (\Omega_N)}
  \lessapprox N^{2\max( (1-\delta)\alpha,
  C\delta)} \|f\|_{L^2} \|g\|_{L^2}
\end{align}
assuming $\operatorname{supp}(\hat{f}) \subset S_1$ and
$\operatorname{supp}(\hat{g}) \subset S_2$ with $\operatorname{diam}(S_j) \le 
1$ and
$\operatorname{dist}(S_1, S_2) \ge c$.


Normalize $f$ and $g$ in $L^2$, and decompose
\[u := U^V_N f = \sum_{T} a_T \phi_T, \quad v := U^V_N = \sum_{T} b_T
  \phi_T\] As in the proof of Lemma~\ref{l:fourier_truncation}, we
discard all but the $O(N^{2d})$ wavepackets whose tubes intersect
$2\Omega_N$.  We also throw away the terms where
$|a_T| = O(N^{-100d})$ or $|b_T| = O(N^{-100d})$, as that portion of
the product can be bounded using the estimates~\eqref{e:wp_bounds} and
Cauchy-Schwartz.

Consequently, in the decompositions of $u$ and $v$ we only consider
the tubes $T$ with $N^{-100d} \lesssim |a_T|, |b_T| \lesssim
1$. Partitioning the interval $[N^{-100d}, 1]$ into $\log N$ dyadic
groups, we may further restrict to the tubes with $|a_T| \sim
\gamma_1$ and $|b_T| \sim \gamma_2$ for dyadic numbers $N^{-100d} \lesssim
\gamma_1, \gamma_2 \lesssim 1$. Let $\mathbf{T}_1$, $\mathbf{T}_2$ be the
tubes for $u$ and $v$, respectively with this property. It therefore
suffices to prove
\begin{align*}
  \Bigl \| \sum_{T_1 \in \mathbf{T}_1} \phi_{T_1} \sum_{T_2 \in \mathbf{T}_2}
  \phi_{T_2} \Bigr\|_{L^{\frac{d+3}{d+1}}(\Omega_N)}
  \lessapprox ( N^{2(1-\delta)\alpha} +
  N^{2C\delta}) \# \mathbf{T}_1^{1/2} \# \mathbf{T}_2^{1/2}
\end{align*}
(we have absorbed the complex phases into the wavepackets).

We have in effect reduced to considering the region of phase space
$\{ (x, \xi) :|x| \lesssim N^2, \ |\xi| \lesssim 1\}$, where the
potential makes only a small perturbation to the Euclidean flow. For
if $|x^s| \lesssim N^2$ and $|t-s| \le N^2$, one has
\begin{align*}
  \begin{split}
  |x^t| &\lesssim N^2\\
  |\xi^t - \xi^s| &\le \int_s^t |\partial_x (V_N)(\tau, x^\tau)| \,
  d\tau \lesssim  \int_s^t |x^\tau| \int_0^1 |\partial^2_xV_N (\tau, s
  x^\tau)| \, ds \, d\tau \lesssim \tau_0 \|
  \partial^2_xV\|_{L^\infty} \lesssim \eta,
\end{split}
\end{align*}
Thus if $\xi \in S_j$, then $\xi^t$ belongs to a small
neighborhood of $S_j$ provided that $\eta \ll c$ is a small multiple
of $c$. For concreteness we choose $\eta$ so that
\begin{align}
  \label{e:freq_constancy}
  |\xi^t - \xi^s|\le \frac{c}{100}.
\end{align}


\subsection{Coarse scale decomposition}

Following Tao, for small $\delta > 0$ we decompose
$\Omega_N = \bigcup_{B \in \mathcal{B}'} B$ into $O(N^{2\delta d})$ smaller
balls of radius $N^{2(1-\delta)}$, and estimate
\begin{align*}
  \Bigl \| \sum_{T_1 \in \mathbf{T}_1} \sum_{T_2 \in \mathbf{T}_2} \phi_{T_1} 
  \phi_{T_2}\|_{L^{\frac{d+3}{d+1}}(\Omega_N)} \lesssim \sum_{B \in 
  \mathcal{B}} 
  \| \sum_{T_1 \in \mathbf{T}_1} \sum_{T_2 \in \mathbf{T}_2} \phi_{T_1} 
  \phi_{T_2}\|_{L^{\frac{d+3}{d+1}}(B)}.
\end{align*}

Let $\sim$ be a relation between tubes and balls to be specified
later. Estimate the norm by the local part
\begin{align}
  \label{e:coarse_local}
\sum_{B \in \mathcal{B}}  \Bigl\| \sum_{T_1 \sim B} \phi_{T_1} \sum_{T_2
  \sim B} \phi_{T_2} \Bigr\|_{L^{\frac{d+3}{d+1}}(B)}
\end{align}
and the global part
\begin{align}
  \label{e:coarse_global}
  \sum_{B \in \mathcal{B}} \Bigl \| \sum_{T_1 \nsim B \text{ or } T_2
  \nsim B} \phi_{T_1} \phi_{T_2} \Bigr\|_{L^{\frac{d+3}{d+1}}(B)}.
\end{align}

We use Proposition~\ref{p:inductive_hyp_cor} with
$R = N^{2(1-\delta)} \le N^2/16$ to estimate the local term by
\begin{align*}
  \eqref{e:coarse_local} &\lessapprox\sum_{B \in \mathcal{B}}
  N^{2(1-\delta)\alpha} \Bigl( \sum_{T_1 \sim B} 1\Bigr)^{1/2} \Bigl(
  \sum_{T_2 \sim B} 1 \Bigr)^{1/2}\\
  &\lessapprox \Bigl( \sum_{T_1 \in \mathbf{T}_1} \# \{B : T_1 \sim B\}
    \Bigr)^{1/2}  \Bigl( \sum_{T_2 \in \mathbf{T}_2} \# \{B : T_2 \sim B\}
    \Bigr)^{1/2}\\
  &\lessapprox 1
\end{align*}
if the relation $\sim$ is chosen so that each $T$ is
associated to $\lessapprox 1$ balls. Note that this step is why we
the Fourier supports are enlarged in
that proposition, as $\operatorname{supp}(\widehat{\phi_{T_1}}(0))$ is not
quite contained in $S_1$.

Heuristically, a judicious choice of $\sim$ allows one to avoid the
worst interactions that would otherwise occur in the bilinear $L^2$
estimate if one were to natively interpolate between $L^1$ and
$L^2$. For example, if all the tubes were to intersect in a
single ball $B$, it would be better to bound
$L^{\frac{d+3}{d+1}}(B)$ directly using the inductive hypothesis rather
than attempt to estimate $L^2(B)$.

The global piece~\eqref{e:coarse_global} is controlled by
interpolating between $L^1$ and $L^2$.  By Cauchy-Schwartz and
conservation of $L^2$ norm,
\begin{align}
  \begin{split}
    \label{e:tube_bilinear_L1}
  &\sum_B \Bigl\| \sum_{T_1 \nsim B \text{ or } T_2 \nsim B}
  \phi_{T_1} \phi_{T_2}
    \Bigr\|_{L^1(B)} \\
  &\lesssim \sum_B  \Bigl( \Bigl \| \sum_{T_1 \sim B} \phi_T
  \Bigr\|_{L^2(B)} +  \Bigl \| \sum_{T_1 \nsim B} \phi_T
  \Bigr\|_{L^2(B)} \Bigr) \Bigl( \Bigl \| \sum_{T_2 \sim B} \phi_T
  \Bigr\|_{L^2(B)} +  \Bigl \| \sum_{T_2 \nsim B} \phi_T
    \Bigr\|_{L^2(B)} \Bigr)\\
    &\lesssim N^{2\delta} N^2 \# \mathbf{T}_1^{1/2} \# \mathbf{T}_2^{1/2}.
    \end{split}
\end{align}

The remaining sections prove the $L^2$ estimate
\begin{align}
  \label{e:tube_bilinear_L2}
  \Bigl \| \sum_{T_1 \nsim B \text{ or } T_2 \nsim B} \phi_{T_1}
  \phi_{T_2} \Bigr\|_{L^2(B)} \lessapprox
  N^{-\frac{d-1}{2}} N^{C\delta}  \# \mathbf{T}_1^{1/2} \mathbf{T}_2^{1/2}.
\end{align}

\subsection{Fine scale decomposition}

Cover $\Omega_N = \bigcup_{q \in \mathbf{q}} q$ by a finitely overlapping 
collection
$\mathbf{q}$ of balls of radius $N$. It suffices to show
\begin{align*}
  \sum_{q \in \mathbf{q}: q \subset 2B}  \Bigl \| \sum_{T_1 \nsim B \text{ or } 
  T_2 \nsim B} \phi_{T_1}
  \phi_{T_2} \Bigr\|_{L^2(q)}^2 \lessapprox
  N^{-(d-1)} N^{C\delta}  \# \mathbf{T}_1 \mathbf{T}_2
\end{align*}

We adopt the following notation from Tao. Fix $q \in \mathbf{q}$ and let
$\mu_1, \mu_2, \lambda_1$ be dyadic numbers.
\begin{itemize}
\item $\mathbf{T}_j(q)$ is the set of tubes $T \in \mathbf{T}_j$ such that $T
  \cap N^\delta q$ is nonempty, where $N^\delta q$ denotes a
  $N^\delta$ neighborhood of $q$.
\item $\mathbf{T}_j^{\nsim B}(q) = \{ T \in \mathbf{T}_j(q) : T \nsim B\}$.
\item $\mathbf{q}(\mu_1, \mu_2)$ is the set of balls $q$ such that
  $ \# \{ T_j \in \mathbf{T}_j : T_j \cap N^\delta q \ne \phi \} \sim
  \mu_j$.

\item $\lambda(T, \mu_1, \mu_2)$ is the number of ($N^\delta$
  neighborhoods of) balls $q \in \mathbf{q}(\mu_1, \mu_2)$ that $T$ intersects.
\item $\mathbf{T}_j[\lambda_1, \mu_1, \mu_2]$ is the set of tubes $T \in
  \mathbf{T}_j$ such that $\lambda(T, \mu_1, \mu_2) \sim \lambda_1$.
\end{itemize}

Pigeonholing dyadically in $\mu_1, \mu_2$, and $\lambda_1$, it
suffices to show
\begin{align*}
  \sum_{q \in \mathbf{q}(\mu_1, \mu_2) : q \subset 2B} \Bigl \| \sum_{T_1
  \in \mathbf{T}_1^{\nsim B}(q) \cap \mathbf{T}_1[\lambda_1, \mu_1, \mu_2]}
  \sum_{T_2 \in \mathbf{T}_2(q)} \phi_{T_1} \phi_{T_2} \Bigr\|_{L^2(q)}^2
  \lessapprox N^{C\delta} N^{-(d-1)} \# \mathbf{T}_1 \# \mathbf{T}_2.
\end{align*}
\subsection{The $L^2$ bound}
\label{ss:restriction_L2_bound}

Fix a ball $q = q(t_q, x_q) \in \mathbf{q}(\mu_1, \mu_2)$ centered at
$(t_q, x_q)$. Suppose want to estimate an expression of the form
\begin{align*}
  \Bigl\| \sum_{T_1} \sum_{T_2} \phi_{T_1} \phi_{T_2} \Bigr\|_{L^2(q)}^2.
\end{align*}

There are two main points to keep in mind:
\begin{itemize}
\item Only tubes that intersect $N^\delta q$ will make a nontrivial
  contribution; that is, tubes whose bicharacteristics $(x^t, \xi^t)$
  satisfy $|x^{t_q} - x_q| \le N^{1+\delta}$.
\item To decouple the contributions of tubes that all overlap near
  $\mathbf{q}$, one needs to exploit oscillation in space and time. While
  Tao employs the spacetime Fourier transform, we instead integrate by
  parts in space and time. 
  Expanding out the $L^2$ norm
  \begin{align}
    \label{e:L2_expanded}
    \sum_{T_1, T_2} \sum_{T_3, T_4} \langle \phi_{T_1} \phi_{T_2}, \,
    \phi_{T_3} \phi_{T_4} \rangle
  \end{align}
  and integrating by parts in both space and time, we shall
  obtain terms of the form
  \[
   (N |\xi_1^{t} + \xi_2^{t} - \xi_3^{t} - \xi_4^{t}|)^{-1}, \quad 
  \bigl(N\bigl| |\xi_1^{t} - \xi_2^{t}|^2 - |\xi_3^{t} - \xi_4^{t}|^2
  \bigr|\bigr)^{-1},
\]
where $(x_j^t, \xi_j^t)$ are bicharacteristics with
$|x_j^{t_q} - x_q| \le N^{1+\delta}$. Since,
by~\eqref{e:integrated_bichar}, the relative frequencies
$\xi_j^t - \xi_k^t$ vary by at most $O(N^{-2+2\delta})$ during the
$O(N^{1+\delta})$ time window when the wavepackets intersect the ball
$N^\delta q$, we can freeze $t = t_q$ above; see
Lemma~\ref{l:kernel_xt_bound} below.
\end{itemize}

Hence, the integral~\eqref{e:L2_expanded} will be small unless
$|x_j^{t_q} - x_q| \le N^{1+\delta}$ for all $j$ \emph{and} the frequencies
$\xi_j^t$ satisfy both resonance conditions
\begin{align}
  \label{e:L2_resonance}
  |\xi_1^{t_q} + \xi_2^{t_q} - \xi_3^{t_q} - \xi_4^{t_q}| = O(N^{-1}), \quad  |\xi_1^{t_q} - \xi_2^{t_q}|^2 - |\xi_3^{t_q} - \xi_4^{t_q}|^2
 = O(N^{-1}).
\end{align}

The preceding discussion motivates the following definition. Let
\begin{align*}
  Z_{q, j} := \{ (x, \xi) : |x| \le 2AN^2, \ \xi \in S_j, \ |x^{t_q} -
  x_q| \le N^{1+\delta} \}.
\end{align*}
For frequencies
$\xi_1$ and $\xi_2'$, define the ``spacetime resonance'' set
\begin{equation}
\nonumber
  \begin{split}
      Z(\xi_1, \xi_2') &= \bigl\{(x_1', \xi_1') \in Z_{q, 1}: \text{ there exists } (x_2, \xi_2) \in Z_{q, 2} \text{ such that }\\
      &\xi_1 + \xi_2^{t_q} = (\xi_1')^{t_q} + \xi_2' \text{ and } |\xi_1 -
      \xi_2^{t_q}|^2 =
      |(\xi_1')^{t_q} -\xi_2'|^2 \bigr\},\\
      \pi(\xi_1, \xi_2') &= \bigl \{ (\xi_1')^{t_q} : (x_1', \xi_1')
      \in Z(\xi_1, \xi_2') \bigr\}.
  \end{split}
\end{equation}
This is a slight modification of Tao's definition which reflects the
time dependence of frequency.

The following lemma follows from elementary geometry.
\begin{lemma}
  \label{l:pi_structure}
  The set $\pi(\xi_1, \xi_2')$ is contained in the hyperplane passing
  through $\xi_1$ and orthogonal to $\xi_2' - \xi_1$ and is therefore
  transverse to $\zeta_2' - \zeta_1$ if $\zeta_1$ and $\zeta_2'$ are
  small perturbations of $\xi_1$ and $\xi_2'$, respectively.
\end{lemma}

Due to the limited time regularity of the phase, we can actually
integrate by parts just once in time. The resulting weaker decay still
turns out to be just enough provided that we slightly refine the
analogue of Tao's main combinatorial estimate for tubes
(estimate~\eqref{e:combinatorial_lma_k} below). Hence we need to account more carefully
for the contributions away from the ``resonant set'' $\pi$.

For $\xi_1, \xi_2'$ and $k > 0$, define the ``time nonresonance'' sets
\begin{align*}
    Z^t_0(\xi_1, \xi_2') &= \bigl\{(x_1', \xi_1')  \in Z_{q, 1}:
                           \text{ there exists } (x_2, \xi_2) \in Z_{q, 2} \text{ such that }\\&|\xi_1 + \xi_2^{t_q} -
                                     (\xi_1')^{t_q} - \xi_2'| \le
  N^{-1+C\delta} \text{ and } \bigl| |\xi_1 - \xi_2^{t_q}|^2 =
                                     |(\xi_1')^{t_q} -\xi_2'|^2\bigr| \le N^{-1+C\delta}\bigr\},\\
    Z^t_k(\xi_1, \xi_2') &= \bigl\{(x_1', \xi_1')  \in Z_{q,1} :
                           \text{ for all } (x_2, \xi_2) \in Z_{q,2}
                           \text{ with } |\xi_1 + \xi_2^{t_q} -
                                     (\xi_1')^{t_q} - \xi_2'| \le
                           N^{-1+C\delta},\\
                                      &\bigl| |\xi_1 - \xi_2^{t_q}|^2 -
                                     |(\xi_1')^{t_q} -\xi_2'|^2\bigr|
  \in (2^{k-1}N^{-1+C\delta},  2^k N^{-1+C\delta}]\bigr\},
\end{align*}
the ``space nonresonance'' set
\begin{align*}
  Z^s(\xi_1, \xi_2') = \bigl\{(x_1', \xi_1')  \in Z_{q,1}:&|\xi_1 + \xi_2^{t_q} -
                                     (\xi_1')^{t_q} + \xi_2'| >
                                                             N^{-1+C\delta}
  \text{ for all } (x_2, \xi_2) \in B_{q,2}\bigr\},
\end{align*}
and the corresponding frequencies at time $t_q$
\begin{align*}
  \pi_k^t(\xi_1, \xi_2') &= \{ (\xi_1')^{t_q} : (x_1', \xi_1') \in
  Z^t_k(\xi_1, \xi_2')\}, \\
  \pi^s(\xi_1, \xi_2') &= \{ (\xi_1')^{t_q} : (x_1', \xi_1') \in
  Z^s(\xi_1, \xi_2')\}.
\end{align*}

An elementary computation shows that
\begin{align}
  \label{e:pi_k_dist}
  \operatorname{dist}(\pi^t_k, \pi) \lesssim 2^k N^{-1+C\delta}.
\end{align}
Indeed, writing $\delta_1 := (\xi_1')^{t_q} - \xi_1$, 
$\delta_2 := \xi_2^{t_q} - \xi_2'$,   and decomposing
$\delta_j = \delta_j^{\parallel} + \delta_j^{\perp}$ into the components
parallel and orthogonal to $\xi_1 - \xi_2'$, we have
\begin{align*}
  | \xi_1 - \xi_2^{t_q}|^2 - | (\xi_1')^{t_q} - \xi_2'|^2 &= |\xi_1 -
                                                            \xi_2' - \delta_2|^2 - |\delta_1 + \xi_1 - \xi_2'|^2\\
                                                          &=
                                                            -2\langle
                                                            \xi_1 -
                                                            \xi_2',
                                                            \delta_1 +
                                                            \delta_2\rangle
                                                            +
                                                            \delta_2^2
                                                            -
                                                            \delta_1^2
  \\
                                                          &= -2 \langle \xi_1 - \xi_2', \delta_1^{\parallel} +
                                                            \delta_2^{\parallel}
                                                            \rangle
                                                            +
                                                            O(N^{-1+C\delta})
                                                            \quad (\text{since \ } |\delta_1-\delta_2| \le N^{1+\delta})\\
                                                          &= -4\langle \xi_1 - \xi_2', \delta_1^\parallel \rangle + O(N^{-1+C\delta}).
\end{align*}
Thus $|
(\xi_1')^{t_q} - \xi_1,\, \xi_1 - \xi_2'\rangle|  \lesssim 2^k
N^{-1+C\delta}$ and the claim follows from
Lemma~\ref{l:pi_structure}.

For $q \in \mathbf{q}(\mu_1, \mu_2)$ with $q \subset 2B$, define
\[
  \mathbf{T}^{\nsim B}_1(q, \lambda_1, \mu_1, \mu_2, \xi_1,
  \xi_2', k)
\]
to be the collection of tubes $T \in \mathbf{T}_1^{\nsim B} (q)
\cap \mathbf{T}_1[\lambda_1, \mu_1, \mu_2]$ such that
$\xi(T)^{t_q} \in \pi_k^t (\xi_1, \xi_2')$. Set
\begin{align}
  \label{e:nuk_def}
  \nu_k(q, \lambda_1, \mu_1, \mu_2) := \sup_{\xi_1 \in S_1, \ \xi_2' \in
  S_2} \# \mathbf{T}_1^{\sim B} (q, \lambda_1, \mu_1, \mu_2, \xi_1^{t_q},
  (\xi_2')^{t_q}, k),
\end{align}
where $|x_1^{t_q} - x_q| + |(x_2')^{t_q} - x_q| \lesssim N^{1+\delta}$.

Then, the analogue of Tao's Lemma 7.1 is:
\begin{lemma}
  \label{l:L2_tube_bound}
  For each $q \in \mathbf{q}(\mu_1, \mu_2)$, we have
  \begin{align*}
    &\Bigl \| \sum_{T_1 \in \mathbf{T}_1^{\nsim B}(q) \cap 
    \mathbf{T}_1[\lambda_1, \mu_1,
              \mu_2]} \sum_{T_2 \in \mathbf{T}_2(q)}
    \phi_{T_1} \phi_{T_2} \Bigr \|_{L^2(q)}^2 \\
    &\lessapprox N^{C\delta}
    N^{-(d-1)} \sup_k 2^{-k} \nu_k(q, \lambda_1, \mu_1, \mu_2) \# 
    (\mathbf{T}_1^{\nsim
    B}(q) \cap \mathbf{T}_1[\lambda_1, \mu_1, \mu_2]) \# \mathbf{T}_2(q).
  \end{align*}
\end{lemma}

\begin{proof}
For conciseness, set
\begin{align*}
  \mathbf{T}_1' &:= \mathbf{T}_1^{\nsim B}(q) \cap \mathbf{T}_1[\lambda_1, 
  \mu_1,
              \mu_2]\\
  \mathbf{T}_2 &:= \mathbf{T}_2(q).
\end{align*}
Then the norm $L^2(q)$ is bounded by the norm $L^2(\eta_N dx dt)$,
where $\eta_N(t)$ is a smooth weight equal to $1$ on $|t-t_q| \le N^{1+\delta}$
and supported in $|t-t_q| \le 2N^{1+\delta}$.
\begin{align*}
  \Bigl\| \sum_{T_1 \in \mathbf{T}_1'} \sum_{T_2 \in \mathbf{T}_2'} \phi_{T_1}
  \phi_{T_2} \Bigr\|_{L^2(\eta_N dx dt)}^2 & = \sum_{T_1, T_1' \in 
  \mathbf{T}_1'}
                               \sum_{T_2, T_2' \in \mathbf{T}_2'} \langle
                               \phi_{T_1} \phi_{T_2} , \phi_{T_1'}
                               \phi_{T_2'}
                               \rangle_{L^2(\chi_N dx dt)}.
\end{align*}
By the bounds~\eqref{e:wp_bounds} and the tranversality of the tubes
in $\mathbf{T}_1'$ and $\mathbf{T}_2'$, the integrand has magnitude
$N^{-2d}$ and is essentially supported on a spacetime ball of width
$N$.  Thus we have the crude bound
\begin{align*}
  | \langle  \phi_{T_1} \phi_{T_2}, \phi_{T_1'} \phi_{T_2'}\rangle|
  \lesssim N^{C\delta} N^{-2d} N^{d+1} = N^{C\delta} N^{-(d-1)}.
\end{align*}
On the other hand, we may integrate by parts to obtain a more refined
bound.
\begin{lemma}
  \label{l:kernel_xt_bound}
  For each $k_1, k_2, \ell \ge 0$ and for all tubes $T_1, T_3 \in
  \mathbf{T}_1'$, $T_2, T_4 \in \mathbf{T}_2'$, we have
\begin{align*}
   | \langle  \phi_{T_1} \phi_{T_2}, \phi_{T_3}
  \phi_{T_4}\rangle|
  \lesssim_{k_1, k_2} N^{C\delta} N^{-(d-1)} &\min \Bigl[ N^{-\ell}|\xi_1^{t_q} +
    \xi_2^{t_q} - \xi_3^{t_q} - \xi_4^{t_q}|^{-\ell},  \\
  &N^{-1} \bigl|
    |\xi_1^{t_q} -\xi_2^{t_q}|^2 - |\xi_3^{t_q} -
    \xi_4^{t_q}|^2 \bigr|^{-1} \Bigr].
\end{align*}
\end{lemma}

\begin{proof}
  The proof has a similar flavor to the earlier
  estimate~\eqref{e:kernel_bound} but takes advantage of oscillation in both space and time.
  
  Let $z_j^t = (x_j^t, \xi_j^t)$ denote the bicharacteristic for
  $\phi_{T_j}$, $j = 1, 2, 3, 4$. By Lemmas~\ref{l:galilei} and
  \ref{l:wp_decomp}, we can write
  \begin{align}
    \label{e:kernel_xt_integral}
    \langle \phi_{T_1} \phi_{T_2}, \phi_{T_3} \phi_{T_4} \rangle =
    \int e^{i\Psi} \phi_1 \phi_2 \overline{\phi_3} \overline{\phi_4}
    \, \eta_N(t) \, dx dt,
  \end{align}
  where $\phi_j$ is a Schr\"{o}dinger solution  which satisfies
  \[ (N \partial_x)^k \phi_j(t, x) \lesssim_{k,M} N^{-d/2} \langle
    N^{-1}(x-x_j^t)\rangle^{-M},\]
  and
  \begin{align*}
    \Psi= \sum_{j=1}^4 \sigma_j \Bigl[ \langle x-x_j^t, \xi_j^t
    \rangle - \int_0^\tau \frac{1}{2}|\xi_j^\tau|^2 - V(\tau, x_j^\tau)
    \, d\tau \Bigr], \quad \sigma = (+, +, -, -).
  \end{align*}

  Using the rapid decay of each $\phi_j$, we may harmlessly (with
  $O(N^{-100d})$ error) localize $\phi_j$ to a $N^{\delta}$
  neighborhood of the tube $T_j$, so that $\phi_j(t)$ is supported in a
  $O(N^{1+\delta})$ neighborhood of the classical path $x_j^t$.
  
  Then
  \begin{align*}
    \partial_x \Psi = \sum_j \sigma_j \xi_j^t, \quad  -\partial_t \Psi = 
    \frac{1}{2} \sum_j \sigma_j|\xi_j^t|^2 + \sum_j \sigma_j
    \bigl[ V(t, x_j^t) + \langle x - x_j^t, \partial_x V(t, x_j^t) \bigr].
  \end{align*}
  The first bound in the statement of the lemma results from
  integrating by parts in $x$, as in the proof
  of~\eqref{e:kernel_bound}, to gain factors of $(N|\xi_1^t + \xi_2^t -
  \xi_3^t - \xi_4^t|)^{-1}$. Since
  \begin{align*}
  \xi_1^t + \xi_2^t -
    \xi_3^t - \xi_4^t =   \xi_1^{t_q} + \xi_2^{t_q} -
    \xi_3^{t_q} - \xi_4^{t_q} + O( N^{-2+2\delta})
  \end{align*}
  during the time window $|t-t_q| \le O(N^{1+\delta})$ when $|x_j^t -
  x_q| \le N^{1+\delta}$, we may replace $t$ by $t_q$.

  As in our work in one space dimension (more specifically,
  the proof of Lemma 4.4 in~\cite{mkv_inv_str_1d}), instead of integrating by
  parts purely in time we use a vector field adapted to the average
  bicharacteristic for the four wavepackets $\phi_{T_j}$. Defining
  \begin{gather*}
    \overline{x^t} : = \frac{1}{4} \sum_{j=1}^4 x_j^t, \quad \overline{\xi^t}
     := \sum_{j=1}^4 \xi_j^t,\\
    L := \partial_t + \langle \overline{\xi^t} , \partial_x \rangle,
  \end{gather*}
  we compute as in that paper that
  \begin{align*}
    -L\Psi &= \frac{1}{2} \sum \sigma_j |\overline{\xi^t}_j|^2 + \sum
    \sigma_j \bigl[ V^{\overline{z}} (t, \overline{x^t}_j) + \langle x-x_j^t,
             \partial_x (V^{\overline{z}})(t, \overline{x^t}_j) \rangle,
  \end{align*}
  where
  \begin{align*}
    \overline{x^t}_j := x_j^t - \overline{x^t}, \quad
    \overline{\xi^t}_j := \xi_j^t - \overline{\xi^t}
  \end{align*}
  denote the coordinates of $\phi_{T_j}(t)$ in phase space relative to
  $(\overline{x^t}, \overline{\xi^t})$; see Figure~\ref{f:fig2}.  
  \begin{figure}
    \includegraphics[scale=0.7]{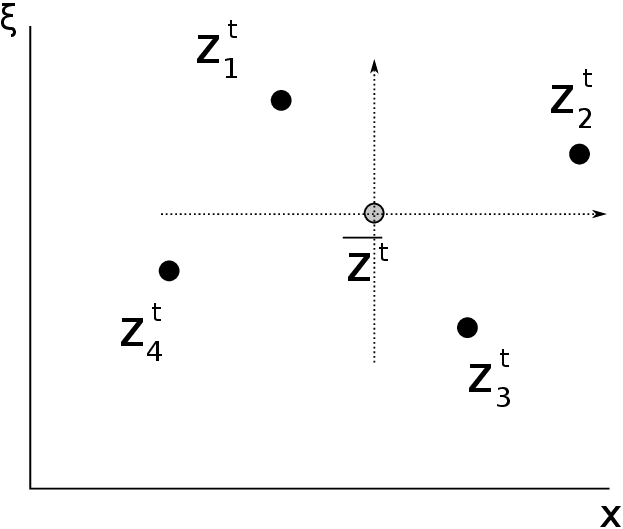}
    \caption{Phase space coordinates relative to the ``center of mass".}
    \label{f:fig2}
  \end{figure}

  We cannot yet integrate by parts since that would require two time
  derivatives of the phase $\Psi$,  but the assumptions on $V$ only allow $\Psi$
  to be differentiated once in time. However, we can decompose
  $\Psi = \Psi_1 + \Psi_2$, where $\Psi_2$ has two time derivatives
  and accounts for the majority of the oscillation of $e^{i\Psi}$;
  indeed, we define $\Psi_1$ and $\Psi_2$ via the ODE
  \begin{align*}
    -L \Psi_2 &= \frac{1}{2} \sum_j \sigma_j |\overline{\xi^t}_j|^2 =  
    \frac{1}{4} \Bigl( |\xi_1^{t_q} - \xi_2^{t_q}|^2 - |\xi_3^{t_q} -
    \xi_4^{t_q}|^2 \Bigr) + O(N^{-2+2\delta}),\\
    -L \Psi_1 &= \sum \sigma_j \Bigl[
    V^{\overline{z}} (t, \overline{x^t}_j) + \langle x - x_j^t,
    \partial_x (V^{\overline{z}}) (t, \overline{x^t}_j) \rangle \Bigr]
    = O(N^{-2+2\delta});
  \end{align*}
  As before we have frozen $t=t_q$ in the main term with
  error at most $O(N^{-2+2\delta})$, and also used the estimates
  $|\overline{x^t}_j|\le \max_{j,k} |x^t_j - x^t_k|\lesssim
  N^{1+\delta}$, $|x-x_j^t| \lesssim N^{1+\delta}$ on the support of
  the integrand~\eqref{e:kernel_xt_integral}.  Note also that the equation
  $\tfrac{d}{dt} \xi_j^t = -\partial_x V(t, x_j^t)$ implies
  $L^2\Psi_2 = O(N^{-2})$. Now integrate by parts using the phase
  $\Psi_2$ to obtain
    \begin{align*}
    \textrm{RH}~\eqref{e:kernel_xt_integral}  &=\int e^{i\Psi_2}
                                                e^{i\Psi_1}\prod_j
                                                \phi_j \, \eta_N(t) \,
                                                dx dt =  i\int e^{i\Psi_2} 
                                                \bigl \langle L, \frac{L 
                                                \Psi_2}{|L\Psi_2|^2}
                                                \rangle e^{i\Psi_1} \phi_1 \phi_2 \overline{\phi_3}
                                                \overline{\phi_4} \,
                                                \eta_N(t) dx dt\\
    &= i \int e^{i\Psi} \Bigl [- \frac{ L^2 \Psi_2}{
      |L\Psi_2|^2}  +  \langle
      \frac{L\Psi_2}{ |L\Psi_2|^2}, i L\Psi_1 + L \rangle\Bigr] \phi_1 \phi_2
      \overline{\phi_3} \overline{\phi_4} \, \eta_N (t) dx dt,
  \end{align*}
  and the second bound in the lemma follows.
\end{proof}

Returning to the proof of Lemma~\ref{l:L2_tube_bound}, we decompose the sum
\begin{align*}
  \sum_{(T_1, T_2') \in \mathbf{T}_1' \times \mathbf{T}_2'} \Bigl[ \sum_{T_1'
  \in \mathbf{T}^s_1} \sum_{T_2 \in \mathbf{T}_2'} + \sum_{0 \le k \lesssim
  \log N} \sum_{T_1' \in \mathbf{T}_{1, k}'} \sum_{T_2 \in \mathbf{T}_2'} 
  \Bigr],
\end{align*}
where $\mathbf{T}^s_{1}$ is the set of tubes in $\mathbf{T}_1'$ whose
bicharacteristic $((x_1')^t, (\xi_1')^t)$ satisfies $(\xi_1')^{t_q} \in
\pi^s ( \xi_1^{t_q}, (\xi_2')^{t_q})$, and we abbreviate
\begin{align*}
  \mathbf{T}_{1,k}' := \mathbf{T}_1^{\nsim B} (q, \lambda_1, \mu_1, \mu_2,
  \xi_1^{t_q}, (\xi_2')^{t_q}, k)
\end{align*}

The contribution from the ``space nonresonance'' terms $\mathbf{T}^s_1$ is
$O(N^{-100d})$.

Now consider the $k$th sum. Lemma~\ref{l:kernel_xt_bound} implies that
\begin{align*}
  | \langle \phi_{T_1} \phi_{T_2}, \, \phi_{T_1'} \phi_{T_2'} \rangle|
  \lesssim N^{C\delta} N^{-(d-1)} 2^{-k}.
\end{align*}
For each $T_1' \in \mathbf{T}_1^{\nsim B}(q, \lambda_1, \mu_1, \mu_2,
\xi_1^{t_q}, (\xi_2')^{t_q}, k)$, the possible tubes $T_2$ correspond
to the bicharacteristics $(x_2^t, x_2^t)$ such that
\begin{align*}
|x_2^{t_q}-x_q| \le N^{1+\delta}, \quad
  \xi_1^{t_q} + \xi_2^{t_q} - (\xi_1')^{t_q} - (\xi_2')^{t_q} =
  O(N^{-1 + C\delta}).
\end{align*}
The preimage of this set under the time $t_q$ Hamiltonian flow map is a
$(N^{1+C\delta})^d \times (N^{-1+C\delta})^{-d}$ box, so there are
$O(N^{C\delta})$ choices of tubes $T_2$. Therefore, the $k$th sum is
at most
\begin{align*}
  N^{C\delta} N^{-(d-1)}2^{-k} \nu_k \# \mathbf{T}_1' \# \mathbf{T}_2',
\end{align*}
whereupon the sum over $k$ is replaced by the supremum at the cost of
a $\log N$ factor.
\end{proof}
It remains to show that
\begin{align}
  \label{e:combinatorial_lma_k}
  \sum_{q \in \mathbf{q}(\mu_1, \mu_2) : q \subset 2B} 2^{-k}\nu_k(q, \lambda_1,
  \mu_1, \mu_2)\# (\mathbf{T}_1^{\nsim B}(q) \cap \mathbf{T}_1[\lambda_1, \mu_1,
              \mu_2]) \# \mathbf{T}_2(q) \lessapprox N^{C\delta} \# \mathbf{T}_1
  \# \mathbf{T}_2.
\end{align}

\subsection{Tube combinatorics}

This section begins exactly as in~\cite[Section 8]{Tao2003}. We define
the relation $\sim$ between tubes and radius $N^{2(1-\delta)}$
balls. For a tube $T \in \mathbf{T}_1[\lambda_1, \mu_1, \mu_2]$, let
$B(T, \lambda_1, \mu_1, \mu_2)$ be a ball $B \in \mathcal{B}$ that
maximizes
\begin{align*}
  \# \{ q \in \mathbf{q}(\mu_1, \mu_2) : T \cap N^\delta q \ne \phi; \ 
  q \cap B \ne \phi\}.
\end{align*}
As $T$ intersects roughly $\lambda_1$ (neighborhoods of) balls
$q \in \mathbf{q}(\mu_1, \mu_2)$ in total and there are $O(N^{2\delta})$
many balls in $\mathcal{B}$, $B(T, \lambda_1, \mu_1, \mu_2)$
must intersect at least $N^{-2\delta} \lambda_1$ of those balls.

Declare $T \sim_{\lambda_1, \mu_1, \mu_2} B'$ if $T \in
\mathbf{T}_1[\lambda_1, \mu_1, \mu_2]$ and $B' \subset 10 B(T, \lambda_1,
\mu_1, \mu_2)$. Finally, for $T \in \mathbf{T}_1$ set $T \sim B$ if $T
\sim_{\lambda_1, \mu_1, \mu_2} B$ for some $\lambda_1, \mu_1,
\mu_2$. Evidently $T \sim B$ for at most $(\log N)^{3} \lessapprox 1 $
many balls.
The relation between tubes in $\mathbf{T}_2$ and
balls in $\mathcal{B}$ is defined similarly.

Now we begin the proof of~\eqref{e:combinatorial_lma_k}. On one hand,
\begin{align*}
  \sum_{q \in \mathbf{q}(\mu_1, \mu_2)} \#( \mathbf{T}_1[\lambda_1, \mu_1,
  \mu_2] \cap \mathbf{T}_1(q))
  &= \sum_{q \in \mathbf{q} (\mu_1, \mu_2)} \sum_{T_1 \in
    \mathbf{T}_1[\lambda_1, \mu_1, \mu_2] \cap \mathbf{T}_1(q)} 1_{T_1 \cap
    N^\delta q \ne 0} \\
  &= \sum_{T \in \mathbf{T}_1[\lambda_1, \mu_1, \mu_2]} \sum_{q \in
    \mathbf{q}(\mu_1, \mu_2)} 1_{T_1 \cap N^\delta q \ne \phi}\\
  &\lesssim \sum_{T \in \mathbf{T}_1} \lambda_1\\
  &= \lambda_1 \# \mathbf{T}_1.
\end{align*}
On the other hand, by definition $\# \mathbf{T}_2(q) \lesssim \mu_2$. The
claim~\eqref{e:combinatorial_lma_k} would therefore follow if we could
show that
\begin{align}
  \label{e:tube_counting_k}
  \nu_k(q_0, \lambda_1, \mu_1, \mu_2) \lessapprox 2^k N^{C\delta} \frac{\#
  \mathbf{T}_2}{\lambda_1 \mu_2}
\end{align}
for all $q_0 \in \mathbf{q}(\mu_1, \mu_2)$ such that $q_0 \subset 2B$.

Fix $\xi_1 \in S_1$, $\xi_2' \in S_2$, and a ball $q_0 = q_0(t_q,
x_q)$. Recalling the definition~\eqref{e:nuk_def} of $\nu_k$, we
need to show that
\begin{align*}
  \# \mathbf{T}^{\nsim B}_1 (q_0, \lambda_1, \mu_1, \mu_2, \xi_1^{t_q},
  (\xi_2')^{t_q}, k) \lessapprox 2^k N^{C\delta} \frac{\# 
  \mathbf{T}_2}{\lambda_1 \mu_2}.
\end{align*}
For brevity write $\mathbf{T}_1' := \mathbf{T}^{\nsim B}_{1,k} (q_0, \lambda_1, 
\mu_1, \mu_2, \xi_1^{t_q},
(\xi_2')^{t_q}, k)$.

Fix $T_1 \in \mathbf{T}_1'$. Since $T_1 \nsim B$, the ball $2B(T_1,
\lambda_1, \mu_1, \mu_2)$ has distance $\gtrsim N^{2(1-\delta)}$ from
$q_0$. Thus
\begin{align*}
  \# \{ q  \in \mathbf{q}(\mu_1, \mu_2) : T_1 \cap N^\delta q \ne \phi, \
  \operatorname{dist}(q, q_0) \gtrapprox N^{2(1-\delta)} \} \gtrapprox
  N^{-2\delta} \lambda_1.
\end{align*}
As each $q \in \mathbf{q}(\mu_1, \mu_2)$ intersects  approximately
$\mu_2$ ($N^\delta$-neighborhoods of) tubes in $\mathbf{T}_2$,
\begin{align*}
    \# \{ (q, T_2)  \in \mathbf{q}(\mu_1, \mu_2) \times \mathbf{T}_2 : T_1
  \cap N^\delta q \ne \phi, \ T_2 \cap N^\delta q \ne \phi, \ 
  \operatorname{dist}(q, q_0) \gtrapprox N^{2(1-\delta)} \} \gtrapprox
  N^{-2\delta} \lambda_1 \mu_2.
\end{align*}
Therefore
\begin{align*}
  &\#\{(q, T_1, T_2) \in \mathbf{q} \times \mathbf{T}_1' \times \mathbf{T}_2 : 
  T_1
  \cap N^\delta q \ne \phi, \ T_2 \cap N^\delta q \ne \phi, \
    \operatorname{dist}(q, q_0) \gtrapprox N^{-2\delta} N^2 \} \\
  &\gtrapprox
  N^{-2\delta} \lambda_1 \mu_2 \#\mathbf{T}_1'
\end{align*}

On the other hand, the cardinality can be bounded above by the
following analogue of Tao's Lemma 8.1:
\begin{lemma}
  \label{l:tube_comb}
  For each  $T_2 \in \mathbf{T}_2$,
  \begin{align*}
    \# \{ (q, T_1) \in \mathbf{q} \times \mathbf{T}_1' : T_1 \cap N^\delta q,
    \ T_2 \cap N^{\delta} q \ne \phi, \ \operatorname{dist}(q, q_0) \gtrapprox
    N^{-2\delta} N^2\} \lessapprox 2^k N^{C\delta}.
  \end{align*}
\end{lemma}

\begin{proof}
  
  We estimate in two steps.

  \begin{itemize}
  \item For any tubes $T_1 \in \mathbf{T}_1'$ and $T_2 \in \mathbf{T}_2$, the
    intersection $N^\delta T_1 \cap N^\delta T_2$ is contained in a
    ball of radius $N^{C\delta}$.
    
  \item The number of tubes $T_1 \in \mathbf{T}_1'$ such that $T_1$
    intersects $N^\delta T_2$ at distance $\gtrapprox N^{-2\delta} N^2$
    from $q_0$ bounded above by $ 2^k N^{C\delta}$.
  \end{itemize}
  The first is evident from transversality. Hence we turn to the
  second claim.

  In Tao's situation, the tubes in $\mathbf{T}_1'$ are all constrained to
  a $O(N^{-1+C\delta})$ neighborhood of a spacetime hyperplane
  transverse to the tube $T_2$ (basically because of
  Lemma~\ref{l:pi_structure}), and there are $O(N^{C\delta})$ many
  such tubes that intersect $T_2$ at distance
  $\gtrapprox N^{-2\delta} N^2$ from $q_0$. The extra $2^k$ factor
  results from the fact that we allow the tubes to deviate from that
  hyperplane by distance $2^k N^{-1+C\delta}$. Also, since our tubes
  are curved it is more convenient to work with their associated
  bicharacteristics instead of using Euclidean geometry in spacetime.

  Fix a tube $T_2 \in \mathbf{T}_2$ with ray
  $t \mapsto (x_2^t, \xi_2^t)$. Then, the tubes $T_1 \in \mathbf{T}_1'$
  such that $N^\delta T_1 \cap N^\delta T_2$ are characterized by the
  property that
  \begin{align*}
    |x(T_1)^t - x_2^t| \lesssim N^{1+\delta} \text{ for some } |t-t_q|
    \gtrapprox N^{-2\delta} N^2.
  \end{align*}
  We need to count the tubes in $\mathbf{T}_1'$ with this property.
  The bicharacteristics for such tubes emanate from the region
  \begin{align*}
    \Sigma := &\{ (x, \xi) : \operatorname{dist}(\xi, S_1) \le N^{-1+C\delta}, \
               \xi^{t_q} \in \pi^t_k, \\
             &|x^{t_q} - x_q| \le N^{1+\delta}, \ |x^t -
               x_2^t| \le N^{1+\delta} \text{ for some } |t-t_q| \gtrsim
               N^{-2\delta} N^2\},
  \end{align*}
  hence it suffices to bound the cardinality of the intersection $(N
  \mathbf{Z}^d \times N^{-1} \mathbf{Z}^d )\cap \Sigma$.
  
  Denote by $\Sigma^{t}$ the image of $\Sigma$ under the time $t$
  Hamiltonian flow map $(x, \xi) \mapsto (x^t, \xi^t)$. Recall
  from~\eqref{e:freq_set_image} that $S_j^t$ denotes the image of the
  initial frequency set $S_j$ for initial positions $x$ with
  $|x| \lesssim N^2$; we saw earlier in~\eqref{e:freq_constancy} that
  $S_j^t$ is a small perturbation of $S_j$.

  Fix a basepoint $x_0$ with $|x_0 - x_q| \le N^{1+\delta}$. By
  Lemma~\ref{l:bichar} and the Hadamard global inverse function
  theorem, when $t \ne t_q$ we can parametrize the graph of the
  flow map $(x^{t_q}, \xi^{t_q}) \mapsto (x^{t}, \xi^{t} )$ by
  the variables
  \[(x^{t_q}, x^{t}) \mapsto \bigl( (x^{t_q}, \xi^{t_q}(x^t, x^{t_q}))
    \mapsto (x^{t}, \xi^{t}( x^{t_q}, x^{t})  )\bigr).\]

  Let $\xi(t,x) := \xi^{t_q}(x_0, x) \in T^*_{x_0}\mathbf{R}^d$ be the
  initial momentum $\xi(t, x) \in T^*_{x_0} \mathbf{R}^d$ such that the
  bicharacteristic with $x^{t_q} = x_0$ and $\xi^{t_q} = \xi(t, x)$
  satisfies $x^t = x$. 
  
  \begin{lemma}
    Suppose at least one $T_1 \in \mathbf{T}_1'$ intersects
    $N^\delta T_2$. For $|t-t_q| \gtrapprox N^{-2\delta} N^2$, the
    curve
    $t \mapsto \zeta_{x_0}(t):= \xi(t, x_2^t) \in T^*_{x_0} \mathbf{R}^d$
    is transverse to the hyperplane containing $\pi(\xi_1, \xi_2')$ for all $\xi_1 \in
    S_1^{t_q}$ and $\xi_2' \in S_2^{t_q}$
    (see Figure~\ref{f:2d_figure}). More precisely there exists
    $C(\eta) > 0$ such that 
    \begin{align*}
      \angle (\dot{\zeta}_{x_0}(t), \pi(\xi_1, \xi_2') ) > C(\eta) \text{ for all } \xi_1
      \in S_1^{t_q}, \ \xi_2' \in S_2^{t_q},
    \end{align*}
    where the angle $\angle (v, W)$ between a vector $v$ and a
    subspace $W$ is defined in the usual manner. Moreover, for each $t$
    the image of a $N^{1+\delta}$ neighborhood of $x_2^t$ under the
    map $x \mapsto \xi(t, x)$ belongs to a $N^{-1+C\delta}$
    neighborhood of $\zeta_{x_0}(t)$.
  \end{lemma}

\begin{figure}
  \includegraphics[scale=0.7]{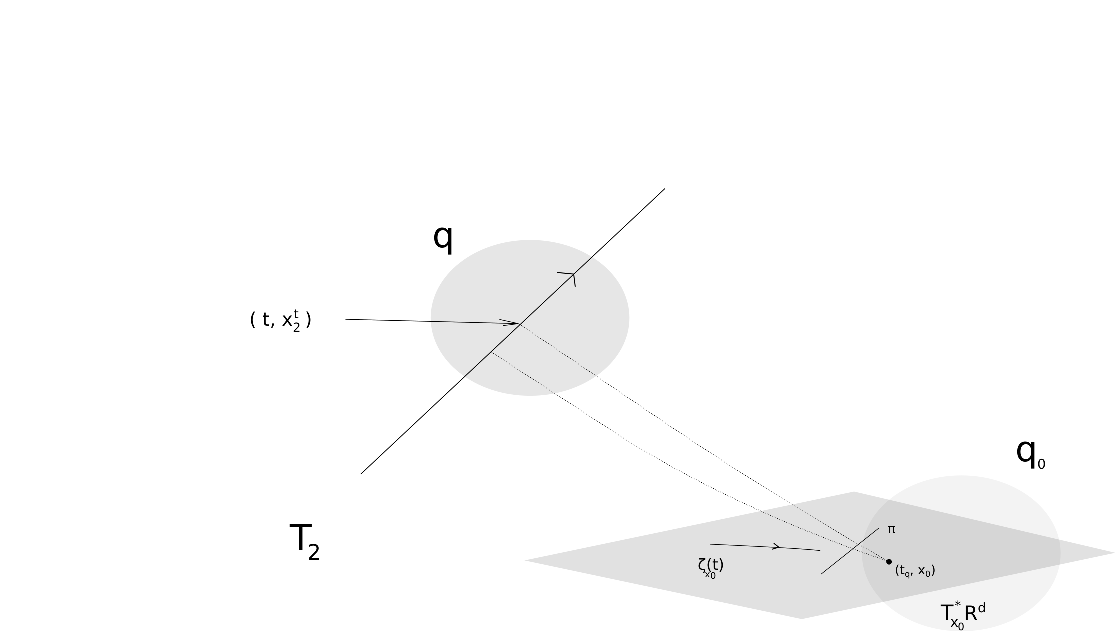}
  \caption{$\zeta_{x_0}(t) \in T^*_{x_0} \mathbf{R}^d$ is the set of
    tangent (covectors) for rays passing through $(t_q, x_0)$ that
    intersect the ray $(t, x_2^t)$ for the tube $T_2$ at times
    $|t-t_q| \gtrapprox N^{2 - 2\delta}$.}
  \label{f:2d_figure}
\end{figure}

\begin{proof}
  By a slight abuse of notation we write
  $\bigl( x^t(y, \zeta), \xi^t(y, \zeta) \bigr)$ for the
  bicharacteristic passing through $(y, \zeta)$ at time $t = t_q$
  instead of $t = 0$.  Both claims are consequences
  of Lemma~\ref{l:bichar}, which yields
  \begin{align*}
    x_2^t &= x^t ( x_0, \zeta_{x_0}(t)), \quad \xi^{t_q} (x_0, \zeta_{x_0}(t)) = \zeta_{x_0}(t),\\
    \xi_2^t &= \frac{d}{dt} x_2^t = \xi^t \bigl( x_0, \zeta_{x_0}(t) \bigr) + 
    \frac{\partial
              x^t}{\partial \zeta_{x_0}} \dot{\zeta}_{x_0}(t)\\
    &= \xi^t(x_0, \zeta_{x_0}(t)) + (t-t_q) \bigl( I + O(\eta) \bigr)
      \dot{\zeta}_{x_0} (t),
  \end{align*}
therefore
\begin{align}
  \label{e:tube_comb_e0}
  \dot{\zeta}_{x_0}(t) = (t-t_q)^{-1} \bigl( I + O(\eta) \bigr) \bigl(
  \xi_2^t - \xi^t(x_0, \zeta_{x_0}(t)) \bigr).
\end{align}

We claim that for any $C > 1$,
\begin{align}
  \label{e:tube_comb_e}
  \operatorname{dist}(\zeta_{x_0} (t_q), S_1^{t_q}) \lesssim_C N^{-1 + C\delta}.
\end{align}
Otherwise, as $|t-t_q| \gtrapprox N^{2(1-\delta)}$, for any ray
$(x_1^t, \xi_1^t)$ with $\xi_1 \in S_1$ and
$|x_1^{t_q} - x_q| \le N^{1+\delta}$, the estimates~\eqref{e:integrated_bichar}
would imply that
\begin{align*}
  |x_1^t - x_2^t| &\gtrsim |t-t_q| |\xi_1^{t_q} - \zeta_{x_0}(t)|  -
                    |x_1^{t_q} - x_0| \\
  &\gtrsim N^{-1+C\delta}  - N^{1+\delta} \gtrsim N^{1 + C\delta},
\end{align*}
so we get the contradiction that every $T_1 \in \mathbf{T}_1'$ misses
$T_2$ by at least $N^{C\delta}$.

By the near-constancy~\eqref{e:freq_constancy} of the frequency
variable and the definition~\eqref{e:small_fourier_supports} of $S_j$,
the covector $\xi_2^t - \xi^t(x_0, \zeta_{x_0}(t))$ belongs a small
perturbation (say, of magnitude at most $\tfrac{c}{50}$) of the
difference set $S_2 - S_1 = -2ce_1 + B(0, \tfrac{c}{50})$, hence by
Lemma~\ref{l:pi_structure} is transverse to the hyperplane
containing $\pi(\xi_1, \xi_2')$. The first claim now follows
from~\eqref{e:tube_comb_e0}.

The argument just given also implies the second statement:  a ray with
$x^{t_q} = x_0$ and $|x_2^t - x^{t}| \le N^{1+\delta}$ must satisfy
$|\xi^{t_q} - \zeta_{x_0}(t)| \lesssim N^{-1+C\delta}$.
\end{proof}

By the second part of the lemma, the fiber of $\Sigma^{t_q}$ in
$T_{x_0}^* \mathbf{R}^d$ is contained in a ``frequency tube''
\begin{align*}
  \Theta(x_0) := \bigcup_{|t-t_q| \gtrapprox N^{2(1-\delta)}} B(\zeta_{x_0}(t), N^{-1+C\delta}).
\end{align*}
As the basepoint $x_0$ varies in $N^{1+\delta}$ neighborhood of $x_q$,
the estimate~\eqref{e:integrated_bichar} implies that the curve
$\zeta_{x_0}(t)$ shifts by at most $O(N^{-1+3\delta})$. Hence the tubes
$\Theta(x_0)$ are all contained in a dilate of $\Theta(x_q)$, which we
denote by
\begin{align*}
  \widetilde{\Theta}(x_q) := \bigcup_{t} B(\zeta_{x_q}(t), N^{-1+C\delta})
\end{align*}
with a larger $C$.

Therefore, $\Sigma^{t_q}$ is contained in the region
\[\tilde{\Sigma}^{t_q} := \bigl\{ (x,\xi) : |x - x_q| \le N^{1+\delta}, \ 
 \xi \in \pi^t_k \cap \widetilde{\Theta} (x_q) \subset \{ \xi \in \widetilde{\Theta}(x_q)
  : \operatorname{dist}(\xi, \pi) \lesssim 2^{k} N^{-1+C\delta}\} \bigr\},
\]
where for the last containment we recall the
estimate~\eqref{e:pi_k_dist}.  The region $\tilde{\Sigma}^{t_q}$ is
sketched in
Figure~\ref{f:2d_figure_2}.
\begin{figure}
  \includegraphics[scale=1.0]{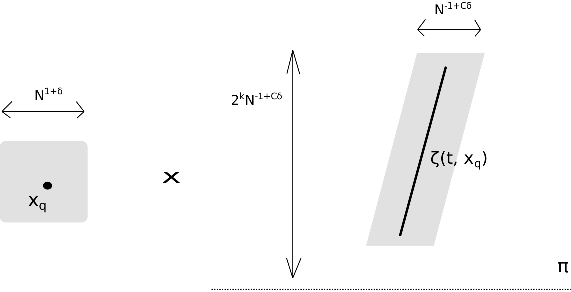}
  \caption{The phase space region $\Sigma^{t_q}$.}
  \label{f:2d_figure_2}
\end{figure}
Using the previous lemma for the central curve $\zeta_{x_q}$, the
frequency projection $(x, \xi) \mapsto \xi)$ of $\tilde{\Sigma}^{t_q}$
can be covered by approximately $2^k$ finitely overlapping cubes
$\bigcup_{1 \le j \lesssim 2^k} Q_j$ of width
$N^{-1+C\delta}$. By~\eqref{e:integrated_bichar}, the preimage of each
box
\begin{align*}
  B(x_q, N^{1+\delta}) \times Q_j
\end{align*}
under the flow map $(x, \xi) \mapsto (x^{t_q}, \xi^{t_q})$ is
contained in a $(CN^{1+C\delta})^d \times (CN^{-1+C\delta})^d$
box. The union of these preimages cover $\Sigma$ and contain at most
$O( 2^k N^{C\delta})$ points in $N \mathbf{Z}^d \times
N^{-1} \mathbf{Z}^{d}$. 
\end{proof}

\section{Remarks on magnetic potentials}
\label{s:magnetic}

We sketch the modifications needed to prove
Theorem~\ref{t:bilinear_Lp_magnetic}. The symbol for $H(t)$ is
\begin{align*}
    a = \frac{1}{2}|\xi|^2 + \langle A , \xi\rangle + V(t, x),
\end{align*}
where $A = A_j(t,x) dx^j$ and $A_j$ are linear functions in the space
variables with bounded time-dependent coefficients.

\begin{itemize}
\item Easy computation shows that the symbol map $a \mapsto a^{z_0}$
  in Lemma~\ref{l:galilei} is
\begin{align*}
  a^{z_0} = \frac{1}{2}|\xi|^2 + \langle A_{(1)}^{z_0}(t,x), \xi \rangle +
  \langle A_{(2)}^{z_0}(t, x), \xi_0^t \rangle + V_{(2)}^{z_0}(t, x),
\end{align*}
where $A_{(1)}^{z_0}(t,x) = A(t, x_0^t + x) - A(t, x_0^t)$ and
$A_{(2)}^{z_0}(t, x) = A(t, x_0^t + x) - \langle x, \partial_x A (t, x_0^t)
\rangle - A(x_0^t)$, and similarly for $V$. Thus when $A$ is linear,
the first order component of the symbol is exactly
``Galilei-invariant'', preserved by the transformation
$a \mapsto a^{z_0}$ in Lemma~\ref{l:galilei}.

\item After
rescaling, the inequality~\eqref{e:bilinear_Lp} takes the form
\begin{align*}
  \| U_N f U_N g\|_{L^{\frac{d+3}{d+1}} ( [-\tau_0 N^2, \tau_0 N^2] \times
  \mathbf{R}^d)} \lesssim_\varepsilon N^\varepsilon \|f\|_{L^2} \|g\|_{L^2},
\end{align*}
where $U_N(t)$ be the propagator for the rescaled symbol
\begin{align*}
  a_N &:= N^{-2} a(N^{-2}t, N^{-1} x, N \xi) = \frac{1}{2}|\xi|^2 +
        N^{-2} \langle A(x), \xi \rangle + N^{-2} V(N^{-2}t, N^{-1}x ).
\end{align*}

\item Exploiting Galilei-invariance, we may reduce to a spatially
  localized estimate as in Proposition~\ref{t:loc_bilinear_Lp}. Note
  that in the region of phase space corresponding to that estimate
  $\{(x, \xi) : |x| \le N^2, \quad |\xi| \lesssim 1\}$, and over a
  $O(N^2)$ time interval, both potential terms have strength $O(1)$
  when integrated over the time interval $|t| \lesssim N^2$. However
  the magnetic term dominates near $x = 0$.

\item Then, the rest of the previous proof can be mimicked with
  essentially no change except for
  Lemma~\ref{l:kernel_xt_bound}. There, one argues essentially as
  before except the vector field $L$ for integrating by parts should
  be replaced by
\begin{align*}
  L := \partial_t + \langle \overline{a_{\xi}(z_j^t)}, \partial_x \rangle,
\end{align*}
where $z_j^t = (x_j^t, \xi_j^t)$ and $\overline{a_{\xi}(z_j^t)} = \tfrac{1}{4} 
\sum_k
a_\xi(z_k^t)$. Then one finds that
\begin{align*}
  -L\Psi = \frac{1}{2} \sum_j \sigma_j |\overline{\xi^t}_j|^2 +
  \sum_j \sigma_j \langle A(\overline{x^t}_j), \xi_j^t \rangle +
  \sum_j \sigma_j \bigl[ V^{\overline{z}} (t, \overline{x^t}_j) +
  \langle x - x_j^t, \partial_x (V^{\overline{z}}) (t,
  \overline{x^t}_j) \rangle \bigr],
\end{align*}
and decomposes as before $\Psi = \Psi_1 + \Psi_2$, where
\begin{align*}
  -L \Psi_1 &= \frac{1}{2} \sum_j \sigma_j |\overline{\xi^t}_j|^2 =
               |\xi_1^{t_q} - \xi_2^{t_q}|^2 - |\xi_3^{t_q} -
              \xi_4^{t_q}|^2 + O(N^{-1+\delta})\\
  -L \Psi_2 &=  \sum_j \sigma_j \langle A(\overline{x^t}_j), \xi_j^t \rangle +
  \sum_j \sigma_j \bigl[ V^{\overline{z}} (t, \overline{x^t}_j) +
  \langle x - x_j^t, \partial_x (V^{\overline{z}}) (t,
              \overline{x^t}_j) \rangle \bigr]\\
  &= O(N^{-1+\delta}).
\end{align*}
\ifdraft
[DRAFT] $|t-s| \le N^{1+\delta}$, $|x_1^s - x_2^s| \le N^{1+\delta}$ imply
\begin{align*}
  |\xi_1^t - \xi_2^t - (\xi_1^s - \xi_2^s)| &\lesssim N^{2+2\delta}
                                              N^{-4} + N^{1+\delta}
                                              N^{-2}\\
                                            &+ N^{4+4\delta} N^{-8} + N^{2+2\delta} N^{-4}\\
  &\lesssim N^{-1+\delta}.
\end{align*}
\fi
As in the proof of Lemma~\ref{l:kernel_xt_bound} the error terms
are computed from the estimates~\eqref{e:integrated_bichar},
$|t-t_q| \lesssim N^{1+\delta}$, and
$|\overline{x^t}_j| \lesssim N^{1+\delta}$. The errors are larger than
before due to the magnetic term $a_{x\xi} = O(N^{-2})$ but are still acceptable.
\end{itemize}

\bibliographystyle{myamsalpha}
\bibliography{bibliography}

\end{document}